\documentclass[12pt]{amsart}
\usepackage{amscd,amssymb,amsthm}

\usepackage[T2A]{fontenc}
\usepackage[cp1251]{inputenc}
\usepackage[english]{babel}

\theoremstyle{plain}
\newtheorem{theorem}{Тheorem}[section]
\newtheorem{corollary}[theorem]{Corollary}
\newtheorem{lemma}[theorem]{Lemma}
\theoremstyle{remark}
\newtheorem{remark}{Remark}[section]
\newtheorem{example}{Example}[section]

\numberwithin{equation}{section}

\renewcommand{\C}{\mathbb{C}}
\newcommand{\N}{\mathbb{N}}
\newcommand{\R}{\mathbb{R}}
\newcommand{\Z}{\mathbb{Z}}

\newcommand{\F}{\mathfrak{F}}

\author[Zastavnyi V.P.]{ Viktor P. Zastavnyi}
\title[Mathieu's series]{Mathieu's series: inequalities, asymptotics and positive definiteness}

\address{Department of Mathematics, Donetsk National University, Universitetskaya str.~24, 340001, Donetsk,  Ukraine}
 \email{zastavn@rambler.ru}


\begin{document}

\begin{abstract}
 Inequalities, asymptotics and, for some specific cases,
 asymptotical expansions were obtained for generalized Mathieu's
 series. A connection between inequalities for Mathieu's series and
 positive definite and completely monotonic functions.
\end{abstract}

\maketitle

\section{Introduction}\label{sec1}

 In 1890   \'{E}mile Leonard Mathieu~\cite{Mathieu} conjectured that
 the inequality
\begin{equation}\label{Mathieu}
S(t):=\sum_{k=1}^{\infty}\frac{2k}{(k^2+t^2)^2}<\frac{1}{t^2} \;,
t>0
\end{equation}
 is valid.
 The first proof of ~\eqref{Mathieu} was published in 1952 in
Berg~\cite{Berg}. It was mentioned that, in case $0<t\le\sqrt{2}$
inequality~\eqref{Mathieu} follows from inequalities
Schr\"oder~\cite[с.~259]{Schroder}
$$
S(t)<\frac{2}{(1+t^2)^2}+\frac{4}{(4+t^2)^2}+\frac{1}{4+t^2}=\frac{1}{t^2}-\frac{2(8-t^6)}{t^2(1+t^2)^2(4+t^2)^2}\;,\;
0<t\le 2\;,
$$
 and, in case  $t\ge\sqrt{2}$ he applied the summation formula of Euler-Maclaurin.

 In 1956 van der Corput and Heflinger~\cite{Corput} got the
 estimation of remainder term in the summation formula of
 Euler-Maclaurin for sufficiently vide classes of the functions. The inequality~\eqref{Mathieu}
 follows from this estimation for every $t>\frac{\sqrt{42}}{16}$.
  If $0<t < 1$, then the function $\frac{x}{(x^2+t^2)^2}$ is strictly decreasing with respect to
   $x\in[1,+\infty)$ and, consequently,
  $S(t)<\frac{2}{(1+t^2)^2}+\int_{1}^{\infty}\frac{2x\;dx}{(x^2+t^2)^2}<\frac{1}{t^2}$.
  They also mentioned the error in
 Emersleben~\cite{Emersleben} which can be corrected. Emersleben has resulted two statements from which
  the inequality~\eqref{Mathieu} follows:
  1) for every natural $t\in\N$ the inequality
   $t^2 S(t)\le 1+\frac{1}{16t^2}$ holds;
 2) the function  $t^2 S(t)$ is strictly increasing for $t>0$.
 The proof of the second statement contains an error.
 The monotonicity of the function  $t^2 S(t)$ was proved in~\cite{Corput}.

In 2007  author obtained the following proof of the
inequality~\eqref{Mathieu}. Let $F(x)=\frac{x}{e^x-1}$, then
$F(0)=1$ and for every $t\ne 0$ the next well-known equality
\begin{equation*}
\begin{split}
  \sum_{k=1}^{\infty}\frac{2k}{(k^2+t^2)^2}=
  &
  \sum_{k=1}^{\infty}\frac{1}{t}\int_{0}^{\infty} e^{-kx}x\,\sin tx\;dx=
  \frac{1}{t}\int_{0}^{\infty} F(x)\,\sin tx\;dx=
  \\&
  \frac{1}{t^2}-\frac{1}{t^2}\int_{0}^{\infty} (-F\,'(x))\,\cos  tx\;dx
\end{split}
\end{equation*}
 valid.
 Using the result of Polya \cite[\S 6, VII]{Polya}, \cite[Theorem 4.3.1]{Lukacs}:
  {\it if $f$
 is convex downwards on $[0,+\infty)$ and $ f(+\infty)=0$, then
$\int_{0}^{\infty} f(x)\cos tx \,dx\ge 0$ for $t\ne 0$}. If, in
addition,  $f$ is not a piecewise linear function with equidistant
nodes then the last integral is strictly greater than zero (see, for
example, Zastavnyi~\cite{Zast2003}). It is easily can be checked
that
 $F\,'(+\infty)=F^{(2)}(+\infty)=0$ and
$F^{(3)}(x)=\frac{e^x}{(e^x-1)^4}\, g(x)$, where
$g(x)=e^{2x}(3-x)-4xe^x-(x+3)=-\sum_{n=5}^{\infty}a_n x^n$,
$a_n=\frac{2^{n-1}(n-6)+4n}{n!}>0$, $ n\ge 5$. Then, for every
 $x>0$ the inequalities $F^{(3)}(x)<0$, $F^{(2)}(x)>0$,
$ F\,'(x)<0$ hold. After application the  Polya's result to the
function  $f(x)=-F\,'(x)$ the inequality~\eqref{Mathieu} follows.

In 1957, Makai~\cite{Makai} proved the inequality~\eqref{Mathieu} by
the following arguments. It is easy to check that, for every
$k\in\N$, $t\ne 0$, the next inequalities
$$
\frac{1}{k(k-1)+t^2+1/2}-\frac{1}{k(k+1)+t^2+1/2}<\frac{2k}{(k^2+t^2)^2}<
\frac{1}{k(k-1)+t^2}-\frac{1}{k(k+1)+t^2}
$$
 hold.
 By the summation, we get
\begin{equation}\label{Mak}
\frac{1}{t^2+a}<\sum_{k=1}^{\infty}\frac{2k}{(k^2+t^2)^2}<\frac{1}{t^2+b}
\;, t>0\;,
\end{equation}
 where $a=\frac{1}{2}$ and $b=0$. The next natural problem appears:
 to find maximal  $b$ and minimal $a$ so that inequality~\eqref{Mak} valid.
  In 1982 Elbert~\cite{Elbert} formulated hypothesis that in  \eqref{Mak}  $a=\frac{1}{2\zeta(3)}$
  can be taken,
 where $\zeta(s)$ is Riemann $\zeta$-function.

In 1998, Alzer, Brenner and Ruehr~\cite{Alzer} proved that
in~\eqref{Mak}  $a=\frac{1}{2\zeta(3)}$ and $b=\frac{1}{6}$ can by
taken and these constants are exact. Their proof is an elementary
consequence of the inequality
\begin{equation}\label{Wilk}
\sum_{k=1}^{\infty}\frac{k}{(k^2+t^2)^3}<
 \left(\sum_{k=1}^{\infty}\frac{k}{(k^2+t^2)^2}\right)^2
 \;,\;t\ge 0\;.
\end{equation}
 The inequality~\eqref{Wilk}, as unsolved problem, was published in 1997~\cite{Alzer97} and proved by Wilkins~\cite{Wilkins}
in 1998. Note that for the calculation of the exact value of
$b=\frac{1}{6}$ the first two terms of the asymptotical expansion
$$
\sum_{k=1}^{\infty}\frac{2k}{(k^2+t^2)^{2}}\sim
\sum_{k=0}^{\infty}\frac{(-1)^kB_{2k}}{t^{2k+2}}\;\;,\;\;t\to+\infty\;,
$$
 was used (\cite{Elbert,Wang}).
 Here, $B_n:=B_n(0)$ - the Bernoulli numbers, and $B_n(x)$ - the Bernoulli
 polynomials.

 We have to mention ones more exact result. In 1980, Diananda
\cite{Diananda} proved the inequality
\begin{equation}\label{Dianand}
 \sum_{k=1}^{\infty}\frac{2k}{(k^2+t^2)^{\mu+1}}<\frac{1}{\mu\, t^{2\mu}} \;,\;\mu>0\;,\; t>0\;,
\end{equation}
 which can be derived from
$$
\frac{2k\mu}{(k^2+t^2)^{\mu+1}}<
\frac{1}{(k(k-1)+t^2)^\mu}-\frac{1}{(k(k+1)+t^2)^\mu}\,,\,\mu>0\,,\,t>0\,,\,k\in\N\,.
$$
 References  on this topic can be found in~\cite{Hoorfar and Qi}.

\section{Main Results}\label{sec2}

 In this paper, we consider the following Mathieu's series
\begin{equation}\label{1}
\begin{split}
 &
 S(t,u,\gamma,\alpha,\mu):=\sum_{k=1}^{\infty}\frac{2(k+u)^{\gamma}}{((k+u)^{\alpha}+t^{\alpha})^{\mu+1}} \;,
 \\&
 \gamma\ge 0\;,\;\alpha>0\;,\;\delta:=\alpha(\mu+1)-\gamma>1\;,\;u>-1 \;,\;t\ge 0\;.
\end{split}
\end{equation}
  Let us introduce the next series
\begin{equation}\label{1a}
\begin{split}
 &
 \tilde{S}(t,u,\gamma,\alpha,\mu):=\sum_{k=1}^{\infty}\frac{2(-1)^{k-1}(k+u)^{\gamma}}{((k+u)^{\alpha}+t^{\alpha})^{\mu+1}} \;,
 \\&
 \gamma\ge 0\;,\;\alpha>0\;,\;\delta:=\alpha(\mu+1)-\gamma>0\;,\;u>-1\;,\;t\ge 0\;.
\end{split}
\end{equation}
 It is obviously, if $\delta>1$,
 \begin{equation}\label{1b}
  \tilde{S}(t,u,\gamma,\alpha,\mu)={S}(t,u,\gamma,\alpha,\mu)-2^{1-\delta}{S}\left(\frac t2,\frac u2,\gamma,\alpha,\mu\right)\;.
   \end{equation}
 In \cite{Qi}, the problem of obtaining the exact inequalities for $S(t,0,\frac{\alpha}{2},\alpha,\mu)$ stated.
 Note that in~\cite[Theorem 2]{Tomovski}, the inequality for  $S(t,0,\frac{\alpha}{2},\alpha,\mu)$
 is not true.

 The Theorem~\ref{thS} (see also Corollaries~\ref{cor1} and \ref{cor2})
 contains the proof of the inequalities for the series~\eqref{1} and~\eqref{1a} for all admissible parameters.
 An asymptotics, as  $t\to+\infty$, and an asymptotic representation, for $(\gamma,\alpha)\in\Z_+\times\N$,
 as a series in the power of $t^{-\alpha(k+\mu+1)}$,
 $k\in\Z_+$, obtained as simple consequences.

 The Bernoulli, $B_n(x)$, and, Euler, $E_n(x)$, polynomials can be
 defined by the following generating functions (see, for example, \cite[\S 1.13, \S 1.14]{Bateman})
 \begin{equation}
\frac{te^{tx}}{e^{t}-1}=\sum_{n=0}^{\infty}\frac{t^n}{n!} \;
B_n(x)\,,\,|t|<2\pi\;; \;\;\;
\frac{2e^{tx}}{e^{t}+1}=\sum_{n=0}^{\infty}\frac{t^n}{n!} \;
E_n(x)\,,\,|t|<\pi\;.
 \end{equation}
 The connection between the Bernoulli and Euler polynomials follows from the next identity
$$
\frac{2e^{tx}}{e^{t}+1}=
 \frac{2e^{tx}}{e^{t}-1}-\frac{4e^{tx}}{e^{2t}-1}
=\sum_{n=1}^{\infty}\frac{t^{n-1}}{(n-1)!} \;\,
\frac{2}{n}\left(B_n(x)-2^nB_n\left(\frac{x}{2}\right)\right)\,,\,0<|t|<\pi\;.
$$
 Consequently,
 $E_{n-1}(x)=\frac{2}{n}\left(B_n(x)-2^nB_n\left(\frac{x}{2}\right)\right)$,
  for every $n\in\N$.
  The Bernoulli and Euler splines defined by the formulas
  $b_n(x)=B_n(\{x\})$ and
  $$e_n(x)=\frac{2}{n+1}\left(b_{n+1}(x)-2^{n+1}b_{n+1}\left(\frac{x}{2}\right)\right)\;,\;n\in\Z_+$$
   correspondingly.

 Theorem~\ref{thEM1} and \ref{thEM2} give the generalization of
 Euler-Maclaurin formula.
\begin{theorem}\label{thEM1}
  {\rm\bf 1.} If, for some $n\in\Z_+$, the function
  $F\in C^n[0,+\infty)\cap L[0,+\infty)$,
  $F^{(n)}$ is a function of bounded variation on $[0,+\infty)$ and
  $F^{(k)}(+\infty)=0$ for $0\le k\le n$, then for every
  $\varepsilon>0$ and $u\ge 0$  the equality
   \begin{equation}\label{EM6}
  \begin{split}
  &\sum_{k=1}^{\infty}F(\varepsilon k+\varepsilon u)=
  \frac{1}{\varepsilon}
  \int_{0}^{\infty} F(t)\,dt +
  \sum_{k=0}^{n}\frac{(-1)^{k}\varepsilon ^k}{(k+1)!}
  \;B_{k+1}(-u)F^{(k)}(0)+
  \frac{(-1)^{n}\varepsilon ^n}{(n+1)!}\;L_n(\varepsilon,u)\;,
  \\&
  \text{ where }
  L_n(\varepsilon,u)=
  \int_{0}^{\infty}b_{n+1}\left(\frac{t}{\varepsilon}\right)\,dF^{(n)}(t+\varepsilon  u) +
  \int_{0}^{1}B_{n+1}(-u+ut)\,dF^{(n)}(t\varepsilon u)
  \;,
  \end{split}
  \end{equation}
    \begin{equation}\label{EM7}
  |L_n(\varepsilon,u)|\le
  \sup_{0\le x\le 1} |b_{n+1}(x)| \;V_{\varepsilon  u}^{\infty}(F^{(n)}) +
  \sup_{-u\le x\le 0} |B_{n+1}(x)| \;V_{0}^{\varepsilon  u}(F^{(n)})
  \end{equation}
  valid.
  If, in addition, $F^{(n)}$ is absolutely continuous on
  $[0,+\infty)$, then $L_n(\varepsilon,u)=o(1)$ for
  $\varepsilon\to+0$,
  uniformly with respect to  $u\in[0,a]$ for every fixed $a>0$.\\
  {\rm\bf 2.} If, for some $q<0$, $n\in\Z_+$, the function
  $F\in C^n[q,+\infty)\cap L[q,+\infty)$,
  $F^{(n)}$ is a function of bounded variation on $[q,+\infty)$ and
  $F^{(k)}(+\infty)=0$, for $0\le k\le n$, then, for every
   $u<0$,  $\varepsilon\in(0,\frac qu)$, the equality~\eqref{EM6}
   valid and
    \begin{equation}\label{EM7a}
  |L_n(\varepsilon,u)|\le
  \sup_{0\le x\le 1} |b_{n+1}(x)| \;V_{\varepsilon  u}^{\infty}(F^{(n)}) +
  \sup_{0\le x\le -u} |B_{n+1}(x)| \;V^{0}_{\varepsilon  u}(F^{(n)})
  \;.
  \end{equation}
  If, in addition, $F^{(n)}$ is absolutely continuous on
  $[q,+\infty)$, then $L_n(\varepsilon,u)=o(1)$, for
  $\varepsilon\to+0$,
  uniformly with respect to  $u\in[b,0]$ for every fixed  $b<0$.
\end{theorem}
\begin{theorem}\label{thEM2}
   {\rm\bf 1.} If, for some  $n\in\Z_+$, the function
  $G\in C^n[0,+\infty)$,
  $G^{(n)}$ is a function of bounded variation on $[0,+\infty)$ and
  $G^{(k)}(+\infty)=0$ for $0\le k\le n$, then, for every
  $\varepsilon>0$ and $u\ge 0$, the equality
   \begin{equation}\label{EM8}
  \begin{split}
  &\sum_{k=1}^{\infty}(-1)^{k-1}G(\varepsilon k+\varepsilon u)=
  \sum_{k=0}^{n}\frac{(-1)^{k}\varepsilon ^k}{2\,k!}
  \;E_{k}(-u)G^{(k)}(0)+
  \frac{(-1)^{n}\varepsilon ^n}{2\,n!}\;l_n(\varepsilon,u)\;,
  \\&
  \text{ where }
  l_n(\varepsilon,u)=
  \int_{0}^{\infty}e_{n}\left(\frac{t}{\varepsilon}\right)\,dG^{(n)}(t+\varepsilon  u) +
  \int_{0}^{1}E_{n}(-u+ut)\,dG^{(n)}(t\varepsilon u)
  \;,
  \end{split}
  \end{equation}
    \begin{equation}\label{EM9}
  |l_n(\varepsilon,u)|\le
  \sup_{0\le x\le 2} |e_{n}(x)| \;V_{\varepsilon  u}^{\infty}(G^{(n)}) +
  \sup_{-u\le x\le 0} |E_{n}(x)| \;V_{0}^{\varepsilon  u}(G^{(n)})
  \end{equation}
  is valid.
  If, in addition,  $G^{(n)}$ is absolutely continuous on
  $[0,+\infty)$, then $l_n(\varepsilon,u)=o(1)$, for
  $\varepsilon\to+0$,
  uniformly on $u\in[0,a]$, for every fixed $a>0$.\\
   {\rm\bf 2.} If, for some $q<0$, $n\in\Z_+$, the function
  $G\in C^n[q,+\infty)$,
  $G^{(n)}$ is a function of bounded variation on $[q,+\infty)$ and
  $G^{(k)}(+\infty)=0$, for $0\le k\le n$, then for every
   $u<0$,  $\varepsilon\in(0,\frac qu)$, the equality~\eqref{EM8} is valid and
    \begin{equation}\label{EM9a}
  |l_n(\varepsilon,u)|\le
  \sup_{0\le x\le 2} |e_{n}(x)| \;V_{\varepsilon  u}^{\infty}(G^{(n)}) +
  \sup_{0\le x\le -u} |E_{n}(x)| \;V^{0}_{\varepsilon  u}(G^{(n)})
  \;.
  \end{equation}
  If, in addition,  $G^{(n)}$ is absolutely continuous on
  $[q,+\infty)$, then $l_n(\varepsilon,u)=o(1)$, for
  $\varepsilon\to+0$, uniformly with respect to $u\in[b,0]$ for every fixed $b<0$.
\end{theorem}
\begin{theorem}\label{thas}
   {\rm\bf 1.} If, for some $q\le 0$, the function
  $F\in C^{\infty}[q,+\infty)$ and $F^{(n)}\in L[q,+\infty)$, for every $n\in\Z_+$,
  then, for every fixed $u\ge 0 $, but if $q<0$, then, for every fixed $u\in\R $,
  the asymptotic representation
   \begin{equation}\label{EM6as}
  \sum_{k=1}^{\infty}F(\varepsilon k+\varepsilon u)\sim
  \frac{1}{\varepsilon}
  \int_{0}^{\infty} F(t)\,dt +
  \sum_{k=0}^{\infty}\frac{(-1)^{k}\varepsilon ^k}{(k+1)!}
  \;B_{k+1}(-u)F^{(k)}(0)\;,\;\varepsilon\to+0\;,
  \end{equation}
  is valid.\\
   {\rm\bf 2.} If, for some $q\le 0$, the function
  $G\in C^{\infty}[q,+\infty)$, $G(+\infty)=0$ and $G^{(n)}\in L[q,+\infty)$, for every $n\in\N$,
  then for every fixed $u\ge 0 $, but if $q<0$, then, for every fixed $u\in\R $,
  the asymptotic representation
 \begin{equation}\label{EM8as}
  \sum_{k=1}^{\infty}(-1)^{k-1}G(\varepsilon k+\varepsilon u)\sim
  \sum_{k=0}^{\infty}\frac{(-1)^{k}\varepsilon ^k}{2\,k!}
  \;E_{k}(-u)G^{(k)}(0)\;,\;\varepsilon\to+0
  \end{equation}
  is valid.
\end{theorem}
\begin{example}\label{ex1}
 The function $F(t)=2t^{2\gamma-1}e^{-t^{2\alpha}}$, for every
 $\gamma,\alpha\in\N$,
 satisfy to the conditions of Theorem~\ref{thas} for every $q<0$. Consequently,
 $F$ has the representation~\eqref{EM6as}. Multiplying~\eqref{EM6as} by $\varepsilon$ and setting
  $\varepsilon^2=x>0$, the next asymptotical representations
 \begin{equation*}
  x^{\gamma}\sum_{k=1}^{\infty}2(k+u)^{2\gamma-1}e^{-(k+u)^{2\alpha} x^{\alpha}}\sim
\frac{\Gamma\left(\frac{\gamma}{\alpha}\right)}{\alpha}+
  \sum_{k=0}^{\infty}\frac{(-1)^{k+1}B_{2(\alpha k+\gamma)}(-u) x^{\alpha k+\gamma}}{(\alpha k+\gamma)k!}
  \;,\;x\to+0\;,
\end{equation*}
\begin{equation*}
  x\sum_{k=1}^{\infty}2(k+u)e^{-(k+u)^2 x}\sim
 \sum_{k=0}^{\infty}\frac{(-1)^{k}B_{2k}(-u) x^k}{k!}
  \;,\;x\to+0
\end{equation*}
 take place, for every fixed $u\in\R $.
    For $u=0$, first of them can be found, for example, in \cite[Example 11.8]{Fedoryuk}.
\end{example}
\begin{theorem}\label{thS}
 Let $\gamma\ge 0$, $\alpha>0$, $\delta:=\alpha(\mu+1)-\gamma>0$ and
 $g(x):=x^{\gamma}(x^\alpha+1)^{-\mu-1}$.
  \begin{enumerate}
                 \item
 If $(\gamma,\alpha)\not\in\Z_+\times\N$, then
 $g\in C^{r}[0,+\infty)$,   $g\not\in C^{r+1}[0,+\infty)$ and, for every integer
 $n\in[0,r]$,
  the function $g^{(n)}$ is absolutely continuous on $[0,+\infty)$,
  where
  \begin{equation}
 r=\left\{
 \begin{array}{ccl}
 [\gamma]&,& \gamma\not\in\Z_+\,,\\
  \gamma+  [\alpha]  &,& \gamma\in\Z_+\,,\, \alpha\not\in\Z_+.
 \end{array}
 \right.
 \end{equation}
 If $\gamma\not\in\Z_+$, then $g^{(p)}(0)=0$ for every integer $p\in[0,r]$.
 If $\gamma\in\Z_+$, $\alpha\not\in\N$, then $g^{(p)}(0)=0$ for every integer
  $p\in[0,r]$, $p\ne\gamma$
 and $g^{(\gamma)}(0)=\gamma!\,$.\\
 Furthermore, the function $G=g$ satisfies to the conditions
  of the statement~{\rm\bf 1} of Theorem~{\rm\ref{thEM2}} for every integer $n\in[0,r]$.\\
 If, in addition,  $\delta>1$, then the function $F=g$ satisfies to the conditions
  of the statement~{\rm\bf 1}
 of Theorem~{\rm\ref{thEM1}} for every integer $n\in[0,r]$.
                 \item
 If $(\gamma,\alpha)\in\Z_+\times\N$, then $g\in C^{\infty}(-1,+\infty)$
 and the function  $G=g$ satisfies to the conditions
  of Theorem~{\rm\ref{thEM2}} for every $q\in(-1,0)$, $n\in\Z_+$, and the function $g^{(n)}$
 is absolutely continuous on $[q,+\infty)$. Furthermore, for every fixed $u\in\R $
 the next asymptotic representation
      \begin{equation}\label{pr1b}
 \tilde{S}(t,u,\gamma,\alpha,\mu)
 \sim
 \sum_{k=0}^{\infty}\frac{(-1)^{k(\alpha+1)+\gamma}}{t^{\alpha(k+\mu+1)}}\cdot
 \frac{\Gamma(\mu+k+1)\,E_{k\alpha+\gamma}(-u)}{\Gamma(\mu+1)\Gamma(k+1)}\;,\;t\to+\infty
\end{equation}
is valid.
 If, in addition,  $\delta>1$, then the function $F=g$ satisfies to the conditions
  of Theorem~{\rm\ref{thEM1}} for every  $q\in(-1,0)$,  $n\in\Z_+$. In this
  case,
  for every fixed $u\in\R$, the next asymptotic representation
      \begin{equation}\label{pr1a}
      \begin{split}
      &
 S(t,u,\gamma,\alpha,\mu)
 \sim
  \frac{1}{t^{\alpha(\mu+1)-\gamma-1}}\cdot
 \frac{2}{\alpha}\,B\left(\frac{\gamma+1}{\alpha},\mu+1-\frac{\gamma+1}{\alpha}\right)+
 \\&
 \sum_{k=0}^{\infty}\frac{(-1)^{k(\alpha+1)+\gamma}}{t^{\alpha(k+\mu+1)}}\cdot
 \frac{2\,\Gamma(\mu+k+1)}{\Gamma(\mu+1)\Gamma(k+1)}\cdot
 \frac{B_{k\alpha+\gamma+1}(-u)}{k\alpha+\gamma+1}\;,\;t\to+\infty\;,
 \end{split}
\end{equation}
 is valid.
  \end{enumerate}
\end{theorem}
 The asymptotic representation   \eqref{pr1b},  for
 $u=0$, $\gamma=1$, $\alpha=2$, and, $\mu=1$, obtained in~\cite{Pogany} by using a different method.

 Further, let us consider some generalization of the inequalities~\eqref{Mak}
 and~\eqref{Dianand}.
 The next Theorem~\ref{thSt} is a sufficient  conditions for the existence of double inequalities.
\begin{theorem}\label{thSt}
 Let the next conditions be satisfied:\\
 {\rm i)} $S(t)\in C[0,+\infty)$, and $S(t)>0$, $t\ge 0$.\\
  {\rm ii)} There exist the positive constants
  $C,\delta_1,A,\beta_1>0$ such that
   \begin{equation}
   S(t)=\frac{C}{t^{\delta_1}}-\frac{A}{t^{\delta_1+\beta_1}}+o\left(\frac{1}{t^{\delta_1+\beta_1}}\right)\;,\;
   t\to+\infty\,.
   \end{equation}

   {\rm iii)}
   For every $t>0$,    $S(t)<\frac{C}{t^{\delta_1}}$.\\
   Then
    \begin{enumerate}
   \item
   If $0<\beta<\beta_1$, the inequality
    $S(t)\le   C(t^\beta+b)^{-\frac{\delta_1}{\beta}}$, $t>0$, is impossible for every $b>0$.
   \item
   If $\beta>\beta_1$, the inequality
    $S(t)\ge C(t^\beta+a)^{-\frac{\delta_1}{\beta}}$, $t>0$, is impossible for every $a\ge0$.
   \item
   Let $m=\inf\limits_{t\ge 0}f(t)$, $M=\sup\limits_{t\ge 0}f(t)$,
   here
    $f(t)=\left(\frac{C}{S(t)}\right)^{\frac{\beta_1}{\delta_1}}-t^{\beta_1}$.
    Then $f(+\infty)=\frac{\beta_1 A}{\delta_1 C}$,
    $0<m\le M<+\infty$, and inequality
     \begin{equation}
  \frac{C}{(t^{\beta_1}+a)^{\frac{\delta_1}{\beta_1}}}
  \le S(t)\le
  \frac{C}{(t^{\beta_1}+b)^{\frac{\delta_1}{\beta_1}}}\;,\; t>0
    \end{equation}
  is valid iff  $0\le b\le m$ and $a\ge M$.
  \end{enumerate}
\end{theorem}
\begin{example}
 For $S(t)=S(t,u,\gamma,\alpha,\mu)$,  $\gamma=0$, $\alpha>0$,
 $\alpha(\mu+1)>1$, $u\ge 0$ the conditions of Theorem~\ref{thSt} be fulfilled.
 In this case  (see Corollary~\ref{cor1})
 $S(t,u,0,\alpha,\mu)=\frac{C}{t^{\delta_1}}-\frac{A}{t^{\delta_1+\beta_1}}+o\left(\frac{1}{t^{\delta_1+\beta_1}}\right)$
 as
 $t\to+\infty$, where
 $C=\frac{2}{\alpha}\,B\left(\frac{1}{\alpha},\mu+1-\frac{1}{\alpha}\right)$,
 $\delta_1=\alpha(\mu+1)-1$, $A=1+2u$, $\beta_1=1$.
 Furthermore, for every $t>0$, the inequality
 \begin{equation*}
 S(t,u,0,\alpha,\mu)<2\int_{0}^{\infty}\frac{dx}{((x+u)^{\alpha}+t^{\alpha})^{\mu+1}}\le
 2\int_{0}^{\infty}\frac{dx}{(x^{\alpha}+t^{\alpha})^{\mu+1}}=
 \frac{C}{t^{\delta_1}}
 \end{equation*}
 is valid.
 Then, for $m=m(u,\alpha,\mu)=\inf\limits_{t\ge 0}f(t)$, and
 $M=M(u,\alpha,\mu)=\sup\limits_{t\ge 0}f(t)$, where
  $$
  f(t)=\left(\frac{\frac{2}{\alpha}\,B\left(\frac{1}{\alpha},\mu+1-\frac{1}{\alpha}\right)}{S(t,u,0,\alpha,\mu)}\right)^{\frac{1}{\alpha(\mu+1)-1}}-t\;,
  $$
  the next inequalities
 \begin{equation*}
 \begin{split}
 &
 \frac{\frac{2}{\alpha}\,B\left(\frac{1}{\alpha},\mu+1-\frac{1}{\alpha}\right)}{(t+M)^{\alpha(\mu+1)-1}}\le
 \sum_{k=1}^{\infty}\frac{2}{((k+u)^{\alpha}+t^{\alpha})^{\mu+1}}
 \le
 \frac{\frac{2}{\alpha}\,B\left(\frac{1}{\alpha},\mu+1-\frac{1}{\alpha}\right)}{(t+m)^{\alpha(\mu+1)-1}}
 \;,\;t>0\;,
 \\&
 0<m\le
   \min\left\{f(0),f(+\infty)\right\}  \le  \max\left\{f(0),f(+\infty)\right\} \le M<+\infty \;,
\end{split}
 \end{equation*}
  is valid. Here
  \begin{equation*}
 f(0)=\left(\frac{\frac{2}{\alpha}\,B\left(\frac{1}{\alpha},\mu+1-\frac{1}{\alpha}\right)}{S(0,u,0,\alpha,\mu)}\right)^{\frac{1}{\alpha(\mu+1)-1}}\;,\;
 f(+\infty)=\frac{\frac 12+u}{(\mu+1-\frac{1}{\alpha}) B\left(\frac{1}{\alpha},\mu+1-\frac{1}{\alpha}\right)}\;.
 \end{equation*}
 In the first inequality  $m$ and $M$ are the best constants. It is easy to see that $m<M$.
\end{example}
 \begin{theorem}\label{cor3}
  Let  $g$ is convex downwards on $(0,+\infty)$ and $\int_{1}^{+\infty}g(x)\,dx$ converges.
  Then, on the open set
 $E:=\{(u,y)\in\R^2:\,(k+u)^2+y>0\,,\,k\in\N\}$,
 the function
 \begin{equation}\label{Suy}
S(u,y):=\sum_{k=1}^{\infty}2(k+u)g((k+u)^2+y)\,,\,(u,y)\in E\,,
 \end{equation}
  is well define and have the next properties:
  \begin{enumerate}
  \item
   $S(u,y)\in C(E)$.
   If $(u_0,y_0)\in E$, then $(u_0,y)\in E$ for every $y>y_0$ and
           $S(u_0,+\infty)=0$.
   For every fixed  $u\ge-1$, the function $S(u,y)$
           is nonnegative, convex downwards with respect to  $y\in(-(1+u)^2,+\infty)$, and decreases.
  If, in addition, $g(x)>0$ for every $x>0$, then, for every fixed
  $u\ge-1$,
  the function $S(u,y)$ is strictly decreasing with respect to
 $y\in(-(1+u)^2,+\infty)$.
  \item
   For every $u\ge -1$, $u^2+u+y>0$, the inequality
\begin{equation}\label{FSF}
F((1+u)^2+y)+\left(\frac{1}{2}+u\right)g((1+u)^2+y)
  \le S(u,y)\le F(u^2+u+y)
\end{equation}
is valid, where
\begin{equation}\label{Ft}
  F(t):=\int_{t}^{+\infty}g(x)\,dx\,,\,t>0\;,\;
  \left(t\ge 0, \text{ if }  \int_{0}^{+\infty}g(x)\,dx<+\infty\right)\,.
 \end{equation}
 Moreover, the right hand side of~\eqref{FSF} going to be the equality iff the function
  $g$ is linear on each segment $[a_k,a_{k+1}]$, $k\in\N$,
 where $a_k=(k-\frac 12+u)^2+y-\frac 14$.
 The left hand side of~\eqref{FSF} valid for every
$u\ge-\frac{3}{2}$, $y>-(1+u)^2$ and going to be the equality only
if the function $g$ is linear on each segment $[s_k,s_{k+1}]$,
$k\in\N$, where $s_k=(k+u)^2+y$.
\end{enumerate}
\end{theorem}
\begin{remark}\label{re1}
Let the function $g$ is convex downwards on $(0,+\infty)$,
$\int_{1}^{+\infty}g(x)\,dx$ converges, and  $g(x)>0$ for every
 $x>0$. It follows from~\eqref{FSF} that, for $u\ge 0$,
\begin{equation}\label{FSF1}
F(a+y)\le S(u,y)\le F(b+y)\;,\;y>0
\end{equation}
 is valid for $a=(1+u)^2$ and $b=u^2+u$. In this case,
 the left hand side of~\eqref{FSF1} is strict inequality.
 If the inequality~\eqref{FSF1} valid for some $a=a_0\ge 0$
  and $b=b_0\ge 0$, then it valid for every  $a\ge a_0$ and $0\le b\le b_0 $, moreover, $b_0\le a_0$.
 Therefore, the problem of definition of the quantities
  \begin{equation}\label{mM}
\begin{split}
 &
 m(u):=\sup\{b\ge 0:\,S(u,y)\le F(b+y)\;,\;y>0\}\,,\,u\ge 0\,,
 \\&
 M(u):=\inf\{a\ge 0:\,F(a+y)\le S(u,y)\;,\;y>0\}\,,\,u\ge 0
\end{split}
  \end{equation}
  naturally appears. It is obviously
  \begin{equation}\label{mM1}
u^2+u\le m(u)\le M(u)\le (1+u)^2\,,\,u\ge 0\,.
  \end{equation}
 Because the function  $F$ strictly decrease, the inverse function $F^{-1}$ exist
 such that it is continuous, positive, and strictly decreasing
 over $I$, where $I=(0,+\infty)$ if $F(+0)=+\infty$ and
$I=(0,F(+0)]$ if $F(+0)<+\infty$. It follows from~\eqref{FSF} that
 \begin{equation}\label{psi}
\psi(u,y):=F^{-1}(S(u,y))-y\,,\,(u,y)\in\R^2_+:=[0,+\infty)\times[0,+\infty)
\end{equation}
 is well defined and, consequently, is continuous on $\R^2_+$. Moreover, the next inequalities
 \begin{equation}\label{225}
u^2+u\le \psi(u,y)<(1+u)^2 \;,\;(u,y)\in\R^2_+
\end{equation}
 is valid. It is obviously, for fixed
$a,b,u\in[0,+\infty)$, the inequality~\eqref{FSF1} is equivalent to
the inequality
 $b\le\psi(u,y)\le a$, $y>0$.
 Therefore,
\begin{equation}\label{mM2}
\begin{split}
 &
 m(u)=\inf_{y\ge 0}\psi(u,y)\,,\,u\ge 0\,,
 \\&
 M(u)=\sup_{y\ge 0}\psi(u,y)\,,\,u\ge 0\,.
\end{split}
  \end{equation}
\end{remark}
\begin{remark}\label{re2}
 Let $g$ is convex downwards on $(0,+\infty)$,
$\int_{1}^{+\infty}g(x)\,dx$ converges, and $g(x)>0$ for every
 $x>0$. Then:

1) $m(u)=u^2+u$ for some $u\ge 0$ $\iff$ $\psi(u,+\infty)=u^2+u$,
 or some  $y\ge 0$ the functions $g$ is linear on each segment $[a_k,a_{k+1}]$, $k\in\N$,
 where $a_k=(k-\frac 12+u)^2+y-\frac 14$.

2) $M(u)=(1+u)^2$ for some $u\ge 0$ $\iff$
$\psi(u,+\infty)=(1+u)^2$.

3) If, for some $u\ge 0$, the equality $m(u)=M(u)=c\ge 0$ is valid,
 then the inequality~\eqref{FSF1} is valid for $a=b=c$. Therefore,
\begin{equation}
m(u)=M(u)=c\ge 0\iff \sum_{k=1}^{\infty}2(k+u)g((k+u)^2+y)\equiv
\int_{c+y}^{+\infty}g(x)\,dx\,,\,y>0\,.
\end{equation}
 For example, if $g(x)=e^{-\lambda x}$, $\lambda>0$, then
 $S(u,0)=\sum_{k=1}^{\infty}2(k+u)e^{-\lambda(k+u)^2}$,
  $$
  S(u,y)=e^{-\lambda y}S(u,0)\;, \;F(t)=\frac{e^{-\lambda  t}}{\lambda}\;,\;
  F^{-1}(s)=-\frac{1}{\lambda}\ln(\lambda s)\;,\;
 \psi(u,y)\equiv -\frac{1}{\lambda}\ln(\lambda S(u,0))
 $$
 and, consequently, $m(u)=M(u)=-\frac{1}{\lambda}\ln(\lambda S(u,0))$.
 In this case the  inequalities~\eqref{225} is strict. Therefore,
  \begin{equation}\label{228}
u^2+u<-\frac{1}{\lambda}\ln\left(\lambda
\sum_{k=1}^{\infty}2(k+u)e^{-\lambda(k+u)^2}\right)<(u+1)^2\;,\;u\ge
0\;,\;\lambda>0\,.
 \end{equation}
\end{remark}
\begin{example}
  If $g(x)=x^{-\mu-1}$, $\mu>0$, $x>0$, then $F(t)=\frac{1}{\mu  t^{\mu}}$,
    $F^{-1}(s)=\left(\frac{1}{\mu s}\right)^{\frac{1}{\mu}}$, $t,s>0$
     and, for every  $u\ge -1$, $u^2+u+y>0$, the inequality
 (see~\eqref{FSF})
\begin{equation}\label{2}
 \frac{1}{\mu((1+u)^2+y)^{\mu}}+
 \frac{\frac{1}{2}+u}{((1+u)^2+y)^{\mu+1}}<
 \sum_{k=1}^{\infty}\frac{2(k+u)}{((k+u)^{2}+y)^{\mu+1}}<\frac{1}{\mu(u^2+u+y)^{\mu}}
\end{equation}
 is valid.
 Left hand side of~\eqref{2} valid for every $u\ge-\frac{3}{2}$, $y>-(1+u)^2$.
\end{example}
 In the case $\gamma=1$, $\alpha=2$, let us denote the series~\eqref{1} by
\begin{equation}\label{Smu}
 S_{\mu}(t,u):=
 S(t,u,1,2,\mu)=\sum_{k=1}^{\infty}\frac{2(k+u)}{((k+u)^2+t^2)^{\mu+1}}
 \;,\;
\mu>0\;,\;u,t\in\R\;.
\end{equation}
  If $-u\in\N$, then, for $k=-u$, the term be absent.
  Therefore, the function $S_{\mu}(t,u)$ is continuous in $t\in\R$ for every fixed $u\in\R$.
  It follows from~\eqref{pr1a} that, for every fixed $u\in\R$, the
  relations
      \begin{equation}\label{asSmu}
      \begin{split}
 &
 S_{\mu}(t,u)
 \sim
 \sum_{k=0}^{\infty}\frac{(-1)^{k}\,B_{2k}(-u)}{t^{2(k+\mu)}}\cdot
 \frac{\Gamma(\mu+k)}{\Gamma(\mu+1)\Gamma(k+1)}
 \;,\;t\to+\infty\;,
 \\&
 S_{\mu}(t,u)=\frac{1}{\mu t^{2\mu}}-
 \frac{u^2+u+\frac 16}{ t^{2\mu+2}}+
 \frac{(u^2+u)^2-\frac{1}{30}}{ t^{2\mu+4}}\cdot\frac{\mu+1}{2}+
 O\left(\frac{1}{ t^{2\mu+6}}\right)
 \;,\;t\to+\infty
 \end{split}
\end{equation}
 hold.
\begin{theorem}\label{thSmu1}
 For every $t>0$, $p>0$, $u>-1$, $\mu>0$ the inequality
 \begin{equation}\label{ner}
 \left|
 \frac{1}{\mu(p^2+t^2)^{\mu}}-S_{\mu}(t,u)
 \right|
 <\left\{
 \begin{array}{ccc}
 \frac{1}{\mu p^{2\mu}}-S_{\mu}(0,u)&,&p-u\le\frac 12\,,\\
  S_{\mu}(0,u)-\frac{1}{\mu p^{2\mu}}&,&p-u\ge 1\,,
 \end{array} \right.
 \end{equation}
 is valid.
 \end{theorem}
\begin{theorem}\label{thSmu2}
 For $\mu>0$, $u\ge 0$, let us set
 $m_{\mu}(u)=\inf\limits_{t\ge 0}f_{\mu,u}(t)$, $M_{\mu}(u)=\sup\limits_{t\ge 0}f_{\mu,u}(t)$,
 where
    $f_{\mu,u}(t)=\left(\mu S_{\mu}(t,u)\right)^{-\frac{1}{\mu}}-t^2$, $t\in\R$.
    Then $f_{\mu,u}(+\infty)=u^2+u+\frac{1}{6}$,
    $0<m_{\mu}(u)\le M_{\mu}(u)<+\infty$ and inequality
     \begin{equation}\label{ner2}
  \frac{1}{\mu(t^2+a)^{\mu}}
  \le S_{\mu}(t,u)\le
  \frac{1}{\mu(t^2+b)^{\mu}}\;,\; t>0
    \end{equation}
  is valid iff  $0\le b\le m_{\mu}(u)$ and $a\ge M_{\mu}(u)$.
  Moreover,  the next statements hold:
  \begin{enumerate}
  \item
  $u^2+u<m_{\mu}(u)\le u^2+u+\frac{1}{6}< u^2+u+\frac{1}{4}<M_{\mu}(u)<(1+u)^2$,
  $\mu>0$, $u\ge 0$.
  \item
  If $\nu>\mu>0$, $u\ge 0$, then
  $m_{\nu}(u)\le m_{\mu}(u)$, $M_{\nu}(u)\ge M_{\mu}(u)$.
  \item
   For every $\nu>0$, $u\ge 0$, the exist $t_{\nu,u}\ge 0$ such that
     $f'_{\nu,u}(t)<0$ for every
   $t>t_{\nu,u}$ and, consequently, the strict inequality
   \begin{equation}\label{ner1}
  \frac{1}{\nu(t^2+f_{\nu,u}(t_{\nu,u}))^{\nu}}
  < S_{\nu}(t,u)<
  \frac{1}{\nu(t^2+f_{\nu,u}(+\infty))^{\nu}}\;,\; t>t_{\nu,u}
    \end{equation}
    is valid.
   \item
   Let us suppose that, for some $\nu>0$, $u\ge 0$, the inequality
   \begin{equation}\label{ner3}
   \left((\nu+1)S_{\nu+1}(t,u)\right)^{\frac{1}{\nu+1}}\le
   \left(\nu S_{\nu}(t,u)\right)^{\frac{1}{\nu}}
   \;,\;t\ge  0
   \end{equation}
   holds.
   Then:\\ {\rm 1)} the function $f_{\nu,u}(t)$ strictly decreases in
   $t\in[0,+\infty)$, the inequality~{\rm\eqref{ner1}} valid for
   $t_{\nu,u}=0$, $M_{\nu}(u)=f_{\nu,u}(0)$ and, for every $\mu\in(0, \nu+1]$,
   the equality  $m_{\mu}(u)=u^2+u+\frac{1}{6}$ and strict
   inequality
  $S_{\mu}(t,u)<\frac{1}{\mu(t^2+m_{\mu}(u))^{\mu}}$, $t\ge 0$ take place.\\
   {\rm 2)}~For every $0<\mu\le\nu$, $0\le a\le u^2+u+\frac{1}{6}$,
  the inequality
  $\frac{d}{dt}\left((t^2+a)^{\mu}S_{\mu}(t,u)\right)>0$, $t>0$ is valid.
  \item
  {\rm 1)} For every $\mu\in(0,2]$,
   the equality $m_{\mu}(0)=\frac{1}{6}$ and strict inequality
   $S_{\mu}(t,0)<\frac{1}{\mu\left(t^2+\frac 16\right)^{\mu}}$, $t\ge 0$, take place.\\
   {\rm 2)}~For every $0<\mu \le 1$, $0\le a\le \frac{1}{6}$,
   $t>0$, the inequality
  $\frac{d}{dt}\left((t^2+a)^{\mu}S_{\mu}(t,0)\right)>0$ is valid.
  \end{enumerate}
\end{theorem}
 As we mention in introduction, the equality $m_{\mu}(0)=\frac{1}{6}$, for
$\mu=1$, proved in~\cite{Alzer}, and the inequality
  $\frac{d}{dt}\left((t^2+a)^{\mu}S_{\mu}(t,0)\right)>0$, $t>0$ for
 $\mu=1$, $a=0$ proved in~\cite{Corput}.

 It follows from Theorem~\ref{thSmu2} that, for $u\ge 0$ the
 relations
\begin{equation}\label{235}
\begin{split}
&
m_\infty(u):=\inf_{\mu>0}m_\mu(u)=\lim_{\mu\to+\infty}m_\mu(u)\;,\;
\\&
M_\infty(u):=\sup_{\mu>0}M_\mu(u)=\lim_{\mu\to+\infty}M_\mu(u)\;,\;
  \\&
  u^2+u\le m_{\infty}(u)\le u^2+u+\frac{1}{6}< u^2+u+\frac{1}{4}<M_{\infty}(u)\le (1+u)^2
\end{split}
\end{equation}
take place.

   The function $f$ is said to be completely monotonic on $(0,+\infty)$ ($f\in M(0,+\infty)$),
   if $f\in  C^{\infty}{(0,+\infty)}$ and for every $k\in\Z_+$, $t>0$,
   the inequality $(-1)^kf^{(k)}(t)\ge 0$ holds.
 The Hausdorff-Bernstein-Widder's Theorem state that $f\in M(0,+\infty)$
 iff
  $f(t)=\int_0^{+\infty}e^{-tx}\ d\mu(x)$, $ t>0$, where $\mu$ is nonnegative Borel measure on
   $[0,+\infty)$ such that the integral converge for every $t>0$.
   Besides, the measure $\mu$ is finite on $[0,+\infty)$ iff $f(+0)<+\infty$.
  \begin{theorem}\label{thSmu3a}
   Let the quantities  $m_{\infty}(u)$, $M_{\infty}(u)$ are defined by~{\rm\eqref{235}}, for $u\ge 0$,
   and
 \begin{equation}\label{236}
 \begin{split}
 &
 \psi_{p,u,\mu}(t):=\frac{1}{\mu(p+t)^{\mu}}-S_{\mu}(\sqrt{t},u)
 \;,\;p\ge 0,\mu>0,t>0\;,\\&
 \varphi_{u}(x):= x\sum_{k=1}^{\infty}2(k+u)e^{-(k+u)^2 x}
 \;,\;x>0\;.
 \end{split}
 \end{equation}
 Then, for each fixed  $u\ge 0$, the next statements hold:
  \begin{enumerate}
  \item
  $\psi_{p,u,\mu}\in M(0,+\infty)$ for every $\mu>0$ $\iff$
    $\psi_{p,u,\mu}\in M(0,+\infty)$ for some $\mu>0$ $\iff$
  $0\le p\le m_{\infty}(u)$.
  \item
  $-\psi_{p,u,\mu}\in M(0,+\infty)$ for every $\mu>0$ $\iff$
    $-\psi_{p,u,\mu}\in M(0,+\infty)$ for some $\mu>0$ $\iff$
  $ p\ge M_{\infty}(u)$.
  \item
   $m_{\infty}(u)=-\sup_{x>0}\frac{\ln\varphi_{u}(x)}{x}$,
   $M_{\infty}(u)=-\inf_{x>0}\frac{\ln\varphi_{u}(x)}{x}$.
    \item
    $u^2+u< m_{\infty}(u)\le u^2+u+\frac{1}{6}$,
    $M_{\infty}(u)=(1+u)^2$.
  \end{enumerate}
  \end{theorem}
  The next corollary immediately follows from Theorem~\ref{thSmu3a}.
\begin{corollary}
   If $u\ge 0$ and $0\le p\le u^2+u$, then
   $\psi_{p,u,\mu}\in M(0,+\infty)$ for every $\mu>0$.
   \end{corollary}
 The function $f:\R^m \to \C$ is said to be positive definite on
$\R^m$, $m\in\N$, if for every $ n\in \N$,  $ \{x_k\}_{k=1}^n
\subset \R^m$ and $ \{c_k\}_{k=1}^n \subset \C$ the next inequality
\begin{equation*}
 \sum_{k,j=1}^n c_k \bar{c}_j f(x_k-x_j) \ge 0
\end{equation*}
 hold.
 About the positive definite functions see, for example, \cite{Zast2000}
 and references therein. For such a function the next inequality
  $|f(x)|\le f(0)$, $x\in\R^m$, holds, and from continuity at zero
  follows continuity in  $\R^m$. By
 Bochner's Theorem the function $f$ is positive definite and continuous in $\R^m$
 iff $f(x)=\int_{\R^m}e^{-i(u,x)} \,d\mu (u),$
 where $\mu$ is nonnegative, finite, Borel's measure on $\R^m$.
 From this Theorem immediately follows the criteria of positive definiteness in
 Fourie's transform term: if
$f\in C(\R^m) \cap L(\R^m),$ then $f$ is positive definite on $\R^m$
iff
\begin{equation*}
\widehat{f}(x)=F_m(f)(x):=\int_{\R^m}e^{i(x,y)}f(y)\ dy\ge 0, \quad
x\in\R^{m}\,,
 \end{equation*}
where $(x,y)=x_1y_1+x_2y_2+...+x_my_m$ - scalar product in $\R^m$,
and, in this case $F_m(f)\in L(\R^m)$.

Let us denote by  $\Phi\left(l_{2}^{m}\right)$ a class of continuous
functions  $f:[0,+\infty)\to\R$  such that function  $f(||x||_2)$ is
positive definite on  $\R^m$. Here $||x||_2=\sqrt{(x,x)}$ the
Euclidean norm in  $\R^m$. Obviously,
$\Phi(l_{2}^{m+1})\subset\Phi(l_{2}^{m})$ and, by well-known
Schoenberg's theorem
\begin{equation}
 \Phi(l_{2}) = \bigcap_{m=1}^{\infty}\Phi(l_{2}^{m}) =\left\{
f(t)=\int_{0}^{+\infty}\exp(-t^{2}s)\ d\mu(s),\ \mu\in
P_{+}\right\}\;,
\end{equation}
where $P_{+}$ - the set of all finite, nonnegative, Borel measures
on $[0,+\infty)$. It follows from the Schoenberg's and
Hausdorff-Bernstein-Widder's theorems that $f\in\Phi(l_{2})$ iff
 $f\in C[0,+\infty)$ and $f(\sqrt{x})\in M(0,+\infty)$.
\begin{theorem}\label{thSmu3}
 Let $m_{\infty}(u)$, $M_{\infty}(u)$ are defined, for $u\ge 0$, by~{\rm\eqref{235}}
 and
 $g_{p,u}(x):=\frac{e^{-p\,x}}{x}-\frac{e^{-ux}}{e^x-1}$
 for $x>0$,
  $g_{p,u}(0):=u+\frac 12-p$, $h_{u}(x):=\frac{e^{-ux}\, x}{e^x-1}$
  for $x>0$, $h_{u}(0):=1$.
  Then, for every fixed  $u\ge 0$ the next statements hold:
  \begin{enumerate}
  \item
  $g_{p,u}\in\Phi(l_{2})\iff 0\le p\le \sqrt{m_{\infty}(u)}$.
  \item
  $-g_{p,u}\in\Phi(l_{2})\iff  p\ge \sqrt{M_{\infty}(u)}=u+1$.
  \item
  The inequality $\frac{d}{dt}\left\{t^{2\mu}S_{\mu}(t,u)\right\}\ge 0$
    is valid for every  $\mu>0$ and $t>0$
   $\iff -\frac{d}{dx}(h_u(x))\in\Phi(l_{2})$.
  \end{enumerate}
\end{theorem}
 The next corollary immediately follows from inequality
   $u^2+u< m_{\infty}(u)$ and Theorem~\ref{thSmu3}.
\begin{corollary}
   If $u\ge 0$ and $0\le p\le \sqrt{u^2+u}$, then
   $g_{p,u}\in\Phi(l_{2})$.
   \end{corollary}
   The proof of theorems of Bochner's, Hausdorff-Bernstein-Widder's and
   Schoenberg's  can be found, for example, in  \cite{TrBel, Vakh, Fell, Scho1, Scho2}.

\section{Euler-Maclaurin's formula. Proofs of Theorems \ref{thEM1}, \ref{thEM2} and \ref{thas}}\label{sec3}

\begin{lemma}\label{lemEMak}
  {\rm\bf 1.} Let $p\in\N$, $u>-p$,  $\varepsilon>0$, $n\in\Z_+$.
  If $F\in C^n[\varepsilon u^{-},\varepsilon(p+u)]$,
  where $u^{-}=\frac{u-|u|}{2}$,  then the next equality
  \begin{equation}\label{EM1}
  \begin{split}
  &\sum_{k=1}^{p}F(\varepsilon k+\varepsilon u)=
  \frac{1}{\varepsilon}
  \int_{0}^{\varepsilon (p+u)} F(t)\,dt +
  \\&
  \sum_{k=0}^{n}\frac{(-1)^{k+1}\varepsilon ^k}{(k+1)!}
  \left(
  B_{k+1}(0)F^{(k)}(\varepsilon p+\varepsilon u)-B_{k+1}(-u)F^{(k)}(0)
  \right)+
  \\&
  \frac{(-1)^{n}\varepsilon ^n}{(n+1)!}
  \left(
  \int_{0}^{\varepsilon p}b_{n+1}\left(\frac{t}{\varepsilon}\right)\,dF^{(n)}(t+\varepsilon  u) +
  \int_{0}^{1}B_{n+1}(-u+ut)\,dF^{(n)}(t\varepsilon u)
  \right)\;
  \end{split}
  \end{equation}
  valid.\\
  {\rm\bf 2.} Let $p\in\N$, $u>-2p$,  $\varepsilon>0$, $n\in\Z_+$.
  If  $G\in C^n[\varepsilon u^{-},\varepsilon(2p+u)]$,
    then the next equality
    \begin{equation}\label{EM2}
  \begin{split}
  &\sum_{k=1}^{2p}(-1)^{k-1}G(\varepsilon k+\varepsilon u)=
  \\&
  \sum_{k=0}^{n}\frac{(-1)^{k+1}\varepsilon ^k}{2\,k!}
  \left(
  E_{k}(0)G^{(k)}(2\varepsilon p+\varepsilon u)-E_{k}(-u)G^{(k)}(0)
  \right)+
  \\&
  \frac{(-1)^{n}\varepsilon ^n}{2\,n!}
  \left(
  \int_{0}^{2\varepsilon p}e_{n}\left(\frac{t}{\varepsilon}\right)\,dG^{(n)}(t+\varepsilon  u) +
  \int_{0}^{1}E_{n}(-u+ut)\,dG^{(n)}(t\varepsilon u)
  \right)\;
  \end{split}
  \end{equation}
  valid.
  Here, $B_n(x)$ and $E_n(x)$ are Bernoulli's and Euler's polynomial, and
  $b_n(x)$ and $e_n(x)$  Bernoulli's and Euler's splines.
\end{lemma}
 \begin{proof}
  Let us consider the statement {\bf 1.}
  Under our conditions, the points $0$ and $\varepsilon u$ belong to the segment $[\varepsilon u^{-},\varepsilon(p+u)]$.
  Consequently, all integrals in~\eqref{EM1} are well defined.
  Obviously, the function $b_n(x)$ are periodic of period $T=1$.
  It follows from the properties of Bernoulli's polynomials that
   $b_0(x)=1$, $b_1(x)=\frac 12+\{x\}$.
    For $n\ge 2$, the functions $b_n(x)$ are absolutely continuous on $\R$,
    and
  $b_n\in C^{n-2}(\R)$,   $b'_n(x)=nb_{n-1}(x)$ for $n\ge 3$,
  $x\in\R$ and $b'_2(x)=2b_1(x)$ for $x\not\in\Z$.
  Taking into account these properties and applying the integration in part in the Riemann-Stieltjes's
  integrals formula, we get
   \begin{equation*}
   \begin{split}
  &
  \frac{(-1)^{n}\varepsilon ^n}{(n+1)!}
  \int_{0}^{\varepsilon p}b_{n+1}\left(\frac{t}{\varepsilon}\right)\,dF^{(n)}(t+\varepsilon  u)=
   \frac{(-1)^{n}\varepsilon ^n}{(n+1)!}b_{n+1}(0)
   \left(
   F^{(n)}(\varepsilon p+\varepsilon u)-F^{(n)}(\varepsilon u)
   \right)
   +\\&
   \frac{(-1)^{n-1}\varepsilon ^{n-1}}{n!}
  \int_{0}^{\varepsilon p}b_{n}\left(\frac{t}{\varepsilon}\right)\,dF^{(n-1)}(t+\varepsilon  u)=
  \ldots=
  \\&
  \sum_{k=0}^{n}\frac{(-1)^{k}\varepsilon ^k}{(k+1)!} b_{k+1}(0)
  \left(
  F^{(k)}(\varepsilon p+\varepsilon u)-F^{(k)}(\varepsilon u)
  \right)-
  \int_{0}^{\varepsilon p} F(t+\varepsilon  u)\,d\,b_{1}\left(\frac{t}{\varepsilon}\right)\;.
  \end{split}
   \end{equation*}
   Taking into account equality  $b_1(x)=-\frac 12+x-[x]$, we
   obtain
\begin{equation*}
   \begin{split}
  \int_{0}^{\varepsilon p} F(t+\varepsilon  u)\,d\,b_{1}\left(\frac{t}{\varepsilon}\right)=
  &
  \frac{1}{\varepsilon}
  \int_{0}^{\varepsilon p} F(t+\varepsilon  u)\,dt-
  \int_{0}^{\varepsilon p} F(t+\varepsilon  u)\,d\,\left[\frac{t}{\varepsilon}\right]=
  \\&
  \frac{1}{\varepsilon}
  \int_{0}^{\varepsilon p} F(t+\varepsilon  u)\,dt-
  \sum_{k=1}^{p}F(\varepsilon k+\varepsilon u)\;.
    \end{split}
   \end{equation*}
   Then
\begin{equation}\label{EM3}
   \begin{split}
  &
  \frac{(-1)^{n}\varepsilon ^n}{(n+1)!}
  \int_{0}^{\varepsilon p}b_{n+1}\left(\frac{t}{\varepsilon}\right)\,dF^{(n)}(t+\varepsilon  u)=
  \\&
  \sum_{k=0}^{n}\frac{(-1)^{k}\varepsilon ^k}{(k+1)!} b_{k+1}(0)
  \left(
  F^{(k)}(\varepsilon p+\varepsilon u)-F^{(k)}(\varepsilon u)
  \right)-
  \frac{1}{\varepsilon}
  \int_{\varepsilon u}^{\varepsilon (p+u)} F(t)\,dt+
  \sum_{k=1}^{p}F(\varepsilon k+\varepsilon u)\;.
  \end{split}
   \end{equation}
    By the same, taking into account the equality  $B'_n(x)=nB_{n-1}(x)$,  $n\in\N$ and
   applying the integration in part formula, we get the next equality
 \begin{equation*}
  \begin{split}
  &
  \frac{(-1)^{n}\varepsilon ^n}{(n+1)!}
   \int_{0}^{1}B_{n+1}(-u+ut)\,dF^{(n)}(t\varepsilon u)=
   \\&
   \sum_{k=0}^{n}\frac{(-1)^{k}\varepsilon ^k}{(k+1)!}
  \left(
  B_{k+1}(0)F^{(k)}(\varepsilon u)-B_{k+1}(-u)F^{(k)}(0)
  \right)-
  \int_{0}^{1} F(t\varepsilon u)\,dB_1(-u+ut)\;.
  \end{split}
  \end{equation*}
   Taking into account equality $B_1(x)=-\frac 12+x$, we obtain
     \begin{equation}\label{EM4}
  \begin{split}
  &
  \frac{(-1)^{n}\varepsilon ^n}{(n+1)!}
   \int_{0}^{1}B_{n+1}(-u+ut)\,dF^{(n)}(t\varepsilon u)=
   \\&
   \sum_{k=0}^{n}\frac{(-1)^{k}\varepsilon ^k}{(k+1)!}
  \left(
  B_{k+1}(0)F^{(k)}(\varepsilon u)-B_{k+1}(-u)F^{(k)}(0)
  \right)-
  \frac{1}{\varepsilon}\int_{0}^{\varepsilon u} F(t)\,dt\;.
  \end{split}
  \end{equation}
 Adding \eqref{EM3} and \eqref{EM4}, taking into account
  $b_k(0)=B_k(0)$, we get equality \eqref{EM1}.

 The statement {\bf 2} follows from the statement {\bf 1} and the
 next equality
 \begin{equation}\label{EM5}
 \sum_{k=1}^{2p}(-1)^{k-1}G(\varepsilon k+\varepsilon u)=
 \sum_{k=1}^{2p}G(\varepsilon k+\varepsilon u)-
 2\sum_{k=1}^{p}G(2\varepsilon k+\varepsilon u)\;
 \end{equation}
 by application the formula  \eqref{EM1}, for $F=G$, to the right
 hand side of \eqref{EM5}:
 to the first term changing $p$ to $2p$,
 to the second term - changing $\varepsilon$ to $2\varepsilon$ and $u$ to $\frac u2$,
 taking into account that $\left(\frac u2\right)^{-}=\frac{^{\;\;}u^{-}}{2}$.
 \end{proof}
\begin{proof}[Proof of Theorem {\rm\ref{thEM1}}]
 At first we will prove the statement {\bf 1.}
 From Lemma~\ref{lemEMak}, by passing to the limit as  $p\to+\infty$
 in \eqref{EM1}, we get \eqref{EM6}.
 From \eqref{EM6} the inequality \eqref{EM7} follows.
 If, in addition,  $F^{(n)}$ is absolutely continuous on $[0,+\infty)$, then
  $F^{(n+1)}\in L[0,+\infty)$ and, consequently,
  $\int_{0}^{\infty}|F^{(n+1)}(t+h)-F^{(n+1)}(t)|\,dt=o(1)$ as $h\to+0$.
  Furthermore, the next equality
  \begin{equation*}
  \begin{split}
  &
  |L_n(\varepsilon,u)|\le
  \left|\int_{0}^{\infty}b_{n+1}\left(\frac{t}{\varepsilon}\right)F^{(n+1)}(t)\,dt\right| +
  \\&
  \sup_{0\le x\le 1} |b_{n+1}(x)|
  \int_{0}^{\infty}\left|F^{(n+1)}(t+\varepsilon u)-F^{(n+1)}(t)\right|\,dt+
  \sup_{-u\le x\le 0} |B_{n+1}(x)| \;V_{0}^{\varepsilon  u}(F^{(n)})
  \;
  \end{split}
  \end{equation*}
  holds.
  Because $b_{n+1}(x)$ is periodic of period $T=1$,
  $F^{(n+1)}\in L[0,+\infty)$, and  $\int_{0}^{1}b_{n+1}(x)\,dx=\int_{0}^{1}\frac{B'_{n+2}(x)}{n+2}\,dx=\frac{B_{n+2}(1)-B_{n+2}(0)}{n+2}=0$
  for every $n\in\Z_+$,
  the Riemann-Lebesgue theorem imply that the first term of last
  inequality tends to zero as $\varepsilon\to+0$.
  The second term tends to zero by arguments above.
   Finally, $V_{0}^{\varepsilon  u}(F^{(n)})\to 0$ as
  $\varepsilon\to+0$ because $F^{(n)}$ is a bounded variation on $[0,+\infty)$
  and continuous at zero function. Thus, $L_n(\varepsilon,u)=o(1)$
  as  $\varepsilon\to+0$, and estimation in $o(1)$ is uniform in  $u\in[0,a]$ for every fixed $a>0$.

  The statement {\bf 2} can be proved in a similar way. Taking in mind that, if in addition
   $F^{(n)}$ is absolutely continuous on
  $[q,+\infty)$, $q<0$, then $F^{(n+1)}\in L[q,+\infty)$ and, consequently,
  $\int_{0}^{\infty}|F^{(n+1)}(t+h)-F^{(n+1)}(t)|\,dt=o(1)$ for
  $h\to-0$ and
  $V^{0}_{\varepsilon  u}(F^{(n)})\to 0$ as
  $\varepsilon\to+0$ since $F^{(n)}$ is a bounded variation on $[q,+\infty)$
  and continuous at zero functions.
\end{proof}
\begin{proof}[Proof of Theorem {\rm\ref{thEM2}}]
  The proof is a similar to the proof of Theorem~\ref{thEM1}.
  Let us prove, for example, the statement {\bf 1.}
  The equality  \eqref{EM8} follows from Lemma~\ref{lemEMak} by passing to the limit in \eqref{EM2}
  as $p\to+\infty$.
  The inequality \eqref{EM9} follows from \eqref{EM8}.
  If, in addition, $G^{(n)}$ is absolutely continuous on
  $[0,+\infty)$, then
  \begin{equation*}
  \begin{split}
  &
  |l_n(\varepsilon,u)|\le
  \left|\int_{0}^{\infty}e_{n}\left(\frac{t}{\varepsilon}\right)G^{(n+1)}(t)\,dt\right| +
  \\&
  \sup_{0\le x\le 2} |e_{n}(x)|
  \int_{0}^{\infty}\left|G^{(n+1)}(t+\varepsilon u)-G^{(n+1)}(t)\right|\,dt+
  \sup_{-u\le x\le 0} |E_{n}(x)| \;V_{0}^{\varepsilon  u}(G^{(n)})
  \;.
  \end{split}
  \end{equation*}
   It is obviously that $e_{n}(x)$ is periodic of period
  $T=2$ function and $\int_{0}^{2}e_{n}(x)\,dx=0$ and so on.
\end{proof}
 \begin{proof}[Proof of Theorem {\rm\ref{thas}}]
  Let us prove the statement {\bf 1.} From the conditions follows the existence of finite limit
  $F^{(n)}(+\infty)=\int_{q}^{\infty}F^{(n+1)}(x)dx+F^{(n)}(q)$ for
  every $n\in\Z_+$.
  Since $F^{(n)}\in L[q,+\infty)$, than $F^{(n)}(+\infty)=0$.
  Further, the Theorem~\ref{thEM1} can be applied to the function
  $F$ for every $n\in\Z_+$.

   The statement {\bf 2} can be proved in a similar way.
 From the conditions follows the existence of finite limit  $G^{(n)}(+\infty)=0$ for
  every $n\in\Z_+$.
  Further, the Theorem~\ref{thEM2} can be applied to the function
  $G$ for every $n\in\Z_+$.
\end{proof}

\section{The inequalities and asymptotic for Mathieu's series. Proof of Theorem \ref{thS} }\label{sec4}

\begin{lemma}\label{lem1}
 Let $\gamma\ge 0$, $\alpha>0$, $\delta:=\alpha(\mu+1)-\gamma>0$,
 $g(x):=x^{\gamma}(x^\alpha+1)^{-\mu-1}$.
 Then
            \begin{enumerate}
 \item
 $g\in C[0,+\infty)\bigcap C^{\infty}(0,+\infty)$
   and $g^{(r)}
   \in L[1,+\infty)$, $g^{(r-1)}(+\infty)=0$ for every $r\in\N$.
   Furthermore, $g\in L[1,+\infty)\iff\delta>1$.
    \item
  If $(\gamma,\alpha)\in\Z_+\times\N$, then $g\in
 C^{\infty}(-1,+\infty)$
 and for every  $p\in\Z_+$ the equality
\begin{equation*}
 \frac{g^{(p)}(0)}{p!}=
 \left\{
 \begin{array}{ccl}
 \frac{(-1)^k\Gamma(\mu+k+1)}{\Gamma(\mu+1)\Gamma(k+1)}
 &,&\text{ если } p=k\alpha+\gamma, k\in\Z_+,\\
 0&,& \text{ если } p\ne k\alpha+\gamma, k\in\Z_+
 \end{array}
 \right.
 \end{equation*}
 valid.
   \item
If $(\gamma,\alpha)\not\in\Z_+\times\N$, then the function $g$ has
finite smoothness of zero
 $$
 r=\left\{
 \begin{array}{ccl}
 [\gamma]&,&\text{ если } \gamma\not\in\Z_+\,,\\
  \gamma+  [\alpha]  &,&\text{ если } \gamma\in\Z_+\,,\, \alpha\not\in\Z_+.
 \end{array}
 \right.
 $$
 In this case, $g\in C^{r}[0,1]$,
  $g\not\in C^{r+1}[0,1]$ and $g^{(r+1)}\in L(0,1)$.
   If $\gamma\not\in\Z_+$ then $g^{(p)}(0)=0$ for every integer $p\in[0,r]$.
 If $\gamma\in\Z_+$, $\alpha\not\in\N$ then $g^{(p)}(0)=0$ for every integer $p\in[0,r]$, $p\ne\gamma$
 and $g^{(\gamma)}(0)=\gamma!$~.
                 \end{enumerate}
\end{lemma}
\begin{proof}
 Proof follows from the next equalities
$$
\frac{1}{(u+1)^{\mu+1}}=\sum_{k=0}^{\infty}\frac{(-1)^k\Gamma(\mu+k+1)}{\Gamma(\mu+1)\Gamma(k+1)}\,u^{
k}\;,\;|u|<1\;,
$$
\begin{equation}\label{lem1a}
\begin{split}
 &
 g(x)=\sum_{k=0}^{\infty}\frac{(-1)^k\Gamma(\mu+k+1)}{\Gamma(\mu+1)\Gamma(k+1)}\cdot x^{\gamma+k\alpha}\;,\;0\le x<1\;,
\\&
 g(x)=\sum_{k=0}^{\infty}\frac{(-1)^k\Gamma(\mu+k+1)}{\Gamma(\mu+1)\Gamma(k+1)}\cdot \frac{1}{x^{\delta+k\alpha}}\;,\; x>1\;.
\end{split}
\end{equation}
\end{proof}
\begin{proof}[Proof of Theorem~{\rm\ref{thS}}]
 follows from Lemma~\ref{lem1},
 the Theorems~\ref{thEM1} and \ref{thEM2} for $G=F=g$ and $\varepsilon=\frac 1t$.
\end{proof}
\begin{corollary}\label{cor1}
 Let $\gamma\ge 0$, $\alpha>0$, $\delta:=\alpha(\mu+1)-\gamma>0$
 and $g(x):=x^{\gamma}(x^\alpha+1)^{-\mu-1}$.
 Then, on the interval $[0,+\infty)$, the $g$ is a bounded variation function which
 equal to $V_0^\infty(g)=1$ if $\gamma=0$ and
  $V_0^\infty(g)=2\left(\frac{\gamma}{\delta}\right)^{\frac{\gamma}{\alpha}}\left(\frac{\gamma}{\delta}+1\right)^{-\mu-1}$
 if $\gamma>0$. Furthermore, the equality
   \begin{equation}
 \tilde{S}(t,u,\gamma,\alpha,\mu)=
 \frac{g(0)+l(t,u)}{t^{\alpha(\mu+1)-\gamma}}\,,\text{ where }
 |l(t,u)|\le V_0^\infty(g)\;,\;l(t,u)=o(1),\,t\to+\infty\,,
   \end{equation}
 holds for every $u\ge 0$ and $t>0$.
 If, in addition, $\delta>1$, then the equality
 \begin{equation}
 S(t,u,\gamma,\alpha,\mu)=
 \frac{\frac{2}{\alpha}\,B\left(\frac{\gamma+1}{\alpha},\mu+1-\frac{\gamma+1}{\alpha}\right)}{t^{\alpha(\mu+1)-\gamma-1}}
  -\frac{(1+2u)\,g(0)-L(t,u)}{t^{\alpha(\mu+1)-\gamma}}
 \end{equation}
 holds for every  $u\ge 0$ and $t>0$, where
 $|L(t,u)|\le (1+2u)V_0^\infty(g)$ and $L(t,u)=o(1)$ as $t\to+\infty$.
 In both cases, the estimation in $o(1)$ is uniform in $u\in[0,a]$ for every fixed $a>0$.
\end{corollary}
\begin{proof}
 In the case $\gamma=0$, the function
  $g$ strictly decrease on $[0,+\infty)$ and, consequently,
 $V_0^\infty(g)=g(0)-g(+\infty)=1$.
 In the case $\gamma>0$, the function
  $g$ increases on $[0,x_0]$ and decreases on $[x_0,+\infty)$,
 where $x_0=\left(\frac{\gamma}{\delta}\right)^{1/\alpha}$ and, consequently,
 $V_0^\infty(g)=2g(x_0)$.
 Further, we apply the statements~{\bf 1} of Theorems~\ref{thEM2} and \ref{thEM1}
 with $n=0$, $F=G=g$, and $\varepsilon=\frac 1t$.
 The inequalities~\eqref{EM8} and \eqref{EM6} with $n=0$ yields
 $l(t,u)=l_0\left(\frac{1}{\varepsilon},u\right)$
 and $L(t,u)=2L_0\left(\frac{1}{\varepsilon},u\right)$.
 Then we use the inequalities~\eqref{EM9} and \eqref{EM7}, having in
 mind that $B_1(x)=-\frac 12 +x$, $b_1(x)=-\frac 12 +\{x\}$,
 $E_0(x)=1$, $e_0(x)=1$ for $0\le x<1$, and $e_0(x)=-1$ for $1\le x<2$,
 and $V_0^{\varepsilon u}(g)+V_{\varepsilon u}^\infty(g)=V_0^\infty(g)$.
\end{proof}
\begin{corollary}\label{cor2}
 Let $(\frac{\gamma}{2},\frac{\alpha}{2})\in\Z_+\times\N$,
 $\delta:=\alpha(\mu+1)-\gamma>0$. Then, for every $n\in\N$ the equality
\begin{equation}
 \tilde{S}(t,0,\gamma,\alpha,\mu)=
 \frac{(-1)^{\gamma}E_{\gamma}(0)}{t^{\alpha(\mu+1)}}+
 o\left(\frac{1}{t^{\alpha(n+\mu+1)}}\right),\,t\to+\infty
   \end{equation}
   holds.
    If, in addition  $\delta>1$, then, for every $n\in\N$, the equality
  \begin{equation}
  \begin{split}
 S(t,0,\gamma,\alpha,\mu)=
 &
 \frac{1}{t^{\alpha(\mu+1)-\gamma-1}}
  \cdot
 \frac{2}{\alpha}\,B\left(\frac{\gamma+1}{\alpha},\mu+1-\frac{\gamma+1}{\alpha}\right)+
 \frac{(-1)^{\gamma}\,2\,B_{\gamma+1}(0)}{(\gamma+1)\,t^{\alpha(\mu+1)}}+
 \\&
  o\left(\frac{1}{t^{\alpha(n+\mu+1)}}\right),\,t\to+\infty
  \end{split}
 \end{equation}
\end{corollary}
 holds.
\begin{proof}
 Proof immediately follows from Theorems~\ref{thEM1}, \ref{thEM2}
 and \ref{thS}, taking into account, that
 $E_m(0)=\frac{2(1-2^{m+1})}{m+1}\,B_{m+1}(0)$, for $m\in\Z_+$,
 $B_1(0)=-\frac 12$, and $B_{2p+1}(0)=0$, for $p\in\N$.
 Thus, $E_0(0)=1$ and $E_{2p}(0)=0$, for  $p\in\N$.
\end{proof}
 \begin{remark}
 It seems to be remainder term in the Corollary~\ref{cor2} approaching zero more foster.
 It can be seen in the next examples.
 With using the summation Poisson's formula, we get
  (see, for example, \cite[Example~11.3]{Fedoryuk}),
 $$
 \sum_{k=1}^{\infty}\frac{2}{k^2+t^2}=\frac{\pi}{t}-\frac{1}{t^2}+\frac{2\pi}{t(e^{2\pi t}-1)}=
 \frac{\pi}{t}-\frac{1}{t^2}+O\left(\frac{1}{te^{2\pi t}}\right)
 \,,\,t\to+\infty\,.
 $$
 From this and the equality \eqref{1b} follows that
$$
 \sum_{k=1}^{\infty}\frac{2(-1)^{k-1}}{k^2+t^2}=\frac{1}{t^2}-\frac{2\pi}{t(e^{\pi t}-e^{-\pi t})}=
 \frac{1}{t^2}+O\left(\frac{1}{te^{\pi t}}\right)
 \,,\,t\to+\infty\,.
 $$
 The last example (see, for example, \cite[\S 3.11--3.12]{Bruijn}):
$$
 \sum_{k=1}^{\infty}\frac{2(-1)^{k-1}}{\sqrt{k^2+t^2}}=\frac{1}{t}+O\left(\frac{1}{\sqrt{t}e^{\pi t}}\right)
 \,,\,t\to+\infty\,.
 $$
 \end{remark}

\section{Hermite-Hadamard's inequality. Proof of Theorems  \ref{thSt} and \ref{cor3}}\label{sec5}

\begin{proof}[Proof of Theorem {\rm\ref{thSt}}]
Let $d\ge 0$, $\beta>0$. Then
 \begin{equation*}
 S(t)\le (\ge) \frac{C}{(t^\beta+d)^{\frac{\delta_1}{\beta}}}\;,\;t>0\iff
 d\le (\ge)
 f_\beta(t):=\left(\frac{C}{S(t)}\right)^{\frac{\beta}{\delta_1}}-t^{\beta}\;,\;t>0\;.
 \end{equation*}
 It is obviously, at $t\to+\infty$, the relations
 \begin{equation*}
 f_\beta(t)=t^{\beta}\left(\left(\frac{C}{S(t)t^{\delta_1}}\right)^{\frac{\beta}{\delta_1}}-1\right)
 \sim
 \frac{t^{\beta}\beta}{\delta_1}\left(\frac{C}{S(t)t^{\delta_1}}-1\right)=
 \frac{t^{\beta}\beta}{\delta_1}\,\frac{C-S(t)t^{\delta_1}}{S(t)t^{\delta_1}}
 \sim
 \frac{t^{\beta-\beta_1}\beta A}{\delta_1 C}
 \end{equation*}
 holds.
 Therefore, if $\beta<\beta_1$, the inequality $0<d\le f_\beta(t)$,
 is impossible for big $t$ (otherwise $0<d\le 0$).
 If $\beta>\beta_1$, the inequality $d\ge f_\beta(t)$
 is impossible for big $t$ (otherwise $d\ge +\infty$).

 Let $\beta=\beta_1$. Then
  $f(t)=f_{\beta_1}(t)>0$ for  $t\ge 0$, and $f(+\infty)=\frac{\beta_1 A}{\delta_1 C}>0$.
  Therefore, $m=\inf\limits_{t\ge 0}f(t)>0$, $M=\sup\limits_{t\ge  0}f(t)<+\infty$.
  Then, for $d\ge 0$, the next statements
 \begin{equation*}
 S(t)\le (\ge) \frac{C}{(t^{\beta_1}+d)^{\frac{\delta_1}{\beta_1}}}\;,\;t>0\iff
 d\le (\ge) f(t)\;,\;t>0\iff d\le m\, (d\ge M)
 \end{equation*}
 are obvious.
\end{proof}
\begin{lemma}\label{ea1}
 Let $g$  is convex downwards on the interval $I$ function, that is, for every
 $x,y\in I$, and $\lambda\in(0,1)$, the inequality $g(\lambda x+(1-\lambda) y)\le \lambda g(x)+(1-\lambda) g(y)$ holds.
  Then, for every $a,b\in I:$ $a<b$, the Hermite-Hadamard's inequality
  \begin{equation}\label{adamar}
 (b-a)\,g\left(\frac{a+b}{2}\right)\le\int_{a}^{b}g(x)\,dx\le(b-a)\,\frac{g(a)+g(b)}{2}
  \end{equation}
  valid.
  The equality in one inequality of \eqref{adamar} occurs iff  $g$ is linear on the segment $[a,b]$.
\end{lemma}
\begin{proof}
 Note, that every convex downwards on the interval $I$ function $f$ is continuous on that interval
 and, for every points
 $x_1,x_2,x_3\in I:$ $x_1<x_2<x_3$ the inequalities
\begin{equation}\label{3hord}
\frac{f(x_2)-f(x_1)}{x_2-x_1}\le
\frac{f(x_3)-f(x_1)}{x_3-x_1}\le\frac{f(x_3)-f(x_2)}{x_3-x_2}
\end{equation}
 hold.
   Therefore, $g\in C(I)$.  The function
$h(x):=g\left(\frac{a+b}{2}+x\right)+g\left(\frac{a+b}{2}-x\right)$
 is convex downwards on $[-c,c]$ with  $c=\frac{b-a}{2}>0$ and, furthermore, is even  function.
 Therefore,  $h$ is increases on $[0,c]$ (it is consequence of \eqref{3hord}
  with $f=h$, $-x_1=x_3\in(0,c]$, and
$x_2\in(0,x_3)$) and, therefore, the inequalities
\begin{equation*}
  \int_{a}^{b}g(x)\,dx=\int_{0}^{c}h(x)\,dx\le c\, h(c)\;(\ge c\, h(0))
\end{equation*}
 hold.
 The inequality  \eqref{adamar} proved. If $g$ is a linear on $[a,b]$,
 then both inequalities of \eqref{adamar} reduces to the equalities and conversely.
 If one inequality of \eqref{adamar} becomes equality then $h$ is a constant on  $[-c,c]$.
  Therefore, the function $g$ is a convex  downwards and upwards  simultaneously on $[a,b]$.
 In this case, for $f=g$, the inequalities in \eqref{3hord} reduces
 to the equalities with $x_1=a$, $x_3=b$ and, consequently,
 $g$ is a linear function on the segment $[a,b]$.
\end{proof}
\begin{lemma}\label{ea2}
 Let $g$  is convex downwards on  $(0,+\infty)$ function such that the integral
  $\int_{1}^{+\infty}g(x)\,dx$ is convergent. Then $g$
 degreases on   $(0,+\infty)$,  $g(+\infty)=0$ and the next statements hold:
   \begin{enumerate}
   \item
 For every $u\ge-1$, $y>-u^2-u$ the inequality
 \begin{equation}\label{adamar1}
\sum_{k=1}^{\infty}2(k+u)\,g((k+u)^2+y)\le
\int_{u^2+u+y}^{+\infty}g(x)\,dx
 \end{equation}
 holds.
 The equality in  \eqref{adamar1} occurs iff the function  $g$ is a linear on each segment
 $[a_{k},a_{k+1}]$, $k\in\N$, with
 $a_{k}=(k-\frac{1}{2}+u)^2+y-\frac{1}{4}$.
   \item
 For every $u\ge-\frac{3}{2}$, $y>-(1+u)^2$ the inequality
 \begin{equation}\label{adamar2}
\int_{(1+u)^2+y}^{+\infty}g(x)\,dx+\left(\frac{1}{2}+u\right)g((1+u)^2+y)
\le\sum_{k=1}^{\infty}2(k+u)\,g((k+u)^2+y)
 \end{equation}
 holds.
 The equality in \eqref{adamar2} occurs iff the function $g$  is a linear on each segment
  $[s_{k},s_{k+1}]$, $k\in\N$, with
 $s_{k}=(k+u)^2+y$.
     \end{enumerate}
\end{lemma}
\begin{proof}
 The inequality \eqref{adamar} implies
\begin{equation*}
\frac{2}{t-a}\,\int^{t}_{a} g(x)\,dx-g(a)\le g(t)\le
\frac{1}{2(t-a)}\,\int^{2t-a}_{a} g(x)\,dx\;,\;t>a>0\,.
\end{equation*}
 Therefore, $\lim_{t\to+\infty}\frac{g(t)}{t}=0$.
 Passing to the limit, as $x_3\to+\infty$, in the inequality~\eqref{3hord} with $f=g$,
 we conclude that $g$ degreases on $(0,+\infty)$. The monotonicity of
$g$ and convergence of the integral $\int_{1}^{+\infty}g(x)\,dx$
imply $g(+\infty)=0$.

 Let us prove the statement {\bf 1}.
 In the left inequality of~\eqref{adamar}, we can take  $a=a_k$ and
 $b=a_{k+1}$, since $a_{k+1}-a_k=2(k+u)\ge 0$ for every $k\in\N$, and $a_1=u^2+u+y>0$.
 Then $\frac{a+b}{2}=(k+u)^2+y$.
 By summation of these inequalities with respect to  $k\in\N$, we get the inequality~\eqref{adamar1}.
 In this case, the equality in~\eqref{adamar1} occurs
(Lemma~\ref{ea1}) only if  $g$ is a linear function on each segment
$[a_{k},a_{k+1}]$, $k\in\N$.

 Let us prove the statement {\bf 2}. In the right inequality of \eqref{adamar},
 we can take $a=s_k$ and $b=s_{k+1}$, since
  $s_{k+1}-s_k=2(k+u)+1\ge 0$,
  for every $k\in\N$, and $s_1=(1+u)^2+y>0$. Then
 \begin{equation*}
 \begin{split}
 \int_{s_k}^{s_{k+1}}g(x)\,dx\le
 &
  \left(k+u+\frac 12\right)
\left(g((k+u)^2+y)+g((k+1+u)^2+y)\right)
 =\\&
(k+u)g((k+u)^2+y)+(k+1+u)g((k+1+u)^2+y)+
 \\&
 \frac 12
(g((k+u)^2+y)-g((k+1+u)^2+y))\,.
\end{split}
\end{equation*}
 We add these inequalities with respect to $k\in\N$.
 By the inequality~\eqref{adamar1}, in the right
 hand side, we get absolutely convergent series
\begin{equation*}
 \int_{s_1}^{+\infty}g(x)\,dx\le
 \sum_{k=1}^{\infty}2(k+u)g((k+u)^2+y)-
 \left(\frac 12+u\right) g((1+u)^2+y)\,.
\end{equation*}
 The last inequality is equivalent to the inequality~\eqref{adamar2},
 and the equality occurs (see Lemma~\ref{ea1}) only if $g$ is a linear function
  on each segment $[s_{k},s_{k+1}]$, $k\in\N$.
\end{proof}
\begin{proof}[Proof of Theorem {\rm\ref{cor3}}]
 Let us prove the property~{\bf 1.} The series \eqref{Suy} converge uniformly on each compact
  $K\subset E$ and, consequently, $S(u,y)\in C(E)$. In fact, for every compact $K\subset E$,
  the number $n_0=n_0(K)$ does exist such that, for every number $n\ge n_0$ and every point $(u,y)\in K$
  the inequalities $n+u\ge-1$ and $(n+u)^2+(n+u)+y\ge n$ hold.
  Hence, the inequality~\eqref{adamar1} can be applied to the remainder of the series~\eqref{Suy}:
   \begin{equation*}
   \begin{split}
   &
\sum_{k=n+1}^{\infty}2(k+u)g((k+u)^2+y)=\sum_{k=1}^{\infty}2(k+n+u)g((k+n+u)^2+y)\le
 \\&
\int\limits^{+\infty}_{(n+u)^2+(n+u)+y}g(x)\,dx\le
\int\limits^{+\infty}_{n}g(x)\,dx\,,\,n\ge n_0\,,\,(u,y)\in K\,.
   \end{split}
   \end{equation*}
 From the last inequality and nonnegativity of  $g$, the uniform convergence of series \eqref{Suy}
 on $K$ follows.

  If
$(u_0,y_0)\in E$ the $(u_0,y)\in E$ for every $y>y_0$ and
$S(u_0,+\infty)=0$.
 Each term of the series  \eqref{Suy} and every their remainder series with the number,
 satisfying  to the inequalities (see above)
 $n+u_0\ge-1$ and $(n+u_0)^2+(n+u_0)+y_0>0$, obeys this property.

 For every fixed $u\ge-1$, the function $S(u,y)$
 is nonnegative, decreasing, and convex downwards with respect to  $y\in(-(1+u)^2,+\infty)$,
 since each term of the series \eqref{Suy} obeys this property.

  If, in addition,  $g(x)>0$ for every $x>0$ then in the inequality \eqref{3hord} for $f=g$,
  we can always choose  $x_3>x_2$ such that the inequality $f(x_3)<f(x_2)$ holds.
 Therefore, the function $g$ strictly decreases on $(0,+\infty)$ and, consequently,
 the function  $S(u,y)$ strictly decreases with respect to
 $y\in(-(1+u)^2,+\infty)$ for every fixed
 $u\ge-1$.

 Note that, if for some  $u\ge-1$ and $y_0>-(1+u)^2$,
 $S(u,y)\equiv const$ for $y\ge y_0$, then $S(u,y)\equiv 0$, for $y\ge y_0$,
 and, consequently, $g(x)\equiv 0$ for $x\ge x_0$, where $x_0=(1+u)^2+y_0$
 if $u>-1$ and $x_0=1+y_0$ if $u=-1$.

 The property~{\bf 2} includes in the Lemma~\ref{ea2}. The Theorem
 \ref{cor3} is proved.
\end{proof}

 \section{$S_{\mu}(t,u)$ series and positive definiteness. Proofs of Theorems  \ref{thSmu1}, \ref{thSmu2}, \ref{thSmu3a} and \ref{thSmu3} }\label{sec6}

\begin{lemma}
 {\rm\bf 1.} For every $(p,t)\in\R^2\setminus\{(0,s):s\le 0\}$, $\nu>\mu>0$, and $u\in\R$,
 the inequality
    \begin{equation}\label{62}
    \begin{split}
    &
  \int_{t}^{+\infty}
  \left(\frac{1}{\nu(p^2+y^2)^{\nu}}-S_{\nu}(y,u)\right)
  y|y^2-t^2|^{\nu-\mu-1}\;dy=
  \\&
  \frac{B(\nu-\mu,\mu+1)}{2}\left(\frac{1}{\mu(p^2+t^2)^{\mu}}-S_{\mu}(t,u)\right)
  \end{split}
  \end{equation}
  holds.\\
     {\rm\bf 2.} For every $u,t,b\in\R$,  $\nu>\mu>0$, the inequality
    \begin{equation}\label{62a}
    \begin{split}
    &
  \int_{t}^{+\infty}
  \left(\nu S_{\nu}(y,u)-(\nu+1)S_{\nu+1}(y,u)(y^2+b)\right)
  y|y^2-t^2|^{\nu-\mu-1}\;dy=
  \\&
  \frac{B(\nu-\mu,\mu+1)}{2}\left(\mu S_{\mu}(t,u)-(\mu+1)S_{\mu+1}(t,u)(t^2+b)\right)
  \end{split}
  \end{equation}
  holds.
\end{lemma}
\begin{proof}
 At first, we prove that if the condition
 \begin{equation*}
  1)\; (p,t)\in\R^2\setminus\{(0,s):s\le 0\}\;,\; \nu>\mu>0
  \; \text{ or }\;
 2)\; p=0\;,\; t<0\;,\; 1>\nu>\mu>0
 \end{equation*}
 holds, then the equality
  \begin{equation}\label{61}
  \int_{t}^{+\infty}
  \frac{y|y^2-t^2|^{\nu-\mu-1}}{(p^2+y^2)^{\nu}}\;dy=
  \frac{B(\nu-\mu,\mu)}{2}\cdot\frac{1}{(p^2+t^2)^{\mu}}
  \end{equation}
  take place.
  It is easy check that the integral in the left hand side of~\eqref{61}
  converges absolutely if the condition 1) or 2) holds.
  In any case, by substitution of variable
  $s=\frac{y^2-t^2}{p^2+t^2}$   in~\eqref{61}, we get
  \begin{equation*}
  \frac{1}{2}\cdot\frac{1}{(p^2+t^2)^{\mu}}
  \int_{0}^{+\infty}\frac{s^{\nu-\mu-1}}{(1+s)^{\nu}}\,ds=
  \frac{B(\nu-\mu,\mu)}{2}\cdot\frac{1}{(p^2+t^2)^{\mu}}\,.
  \end{equation*}
  The equality~\eqref{61} proved. After substitution in~\eqref{61} $\nu$
  by  $\nu+1$ and  $\mu$ by $\mu+1$, the equality
   \begin{equation*}
  \int_{t}^{+\infty}
  \frac{y|y^2-t^2|^{\nu-\mu-1}}{(p^2+y^2)^{\nu+1}}\;dy=
   \frac{B(\nu-\mu,\mu+1)}{2}\cdot\frac{1}{(p^2+t^2)^{\mu+1}}
  \end{equation*}
  take place for every
  $(p,t)\in\R^2\setminus\{(0,s):s\le 0\}$,   $\nu>\mu>-1$.
  Tacking into account the last equality, for every $u, t\in\R$,
  $\nu>\mu>0$, we get
    \begin{equation}\label{61a}
    \begin{split}
    &
  \int_{t}^{+\infty}
  S_{\nu}(y,u)
  y|y^2-t^2|^{\nu-\mu-1}\;dy=
  \\&
  \sum_{k=1}^{\infty}2(k+u)\int_{t}^{+\infty}
  \frac{y|y^2-t^2|^{\nu-\mu-1}}{((k+u)^2+y^2)^{\nu+1}}\;dy=
  \frac{B(\nu-\mu,\mu+1)}{2}S_{\mu}(t,u)\,.
  \end{split}
  \end{equation}
  Note that, in the series \eqref{Smu} for $S_{\nu}(y,u)$ and
  $S_{\mu}(t,u)$, the terms for $k=-u$ we absent if $-u\in\N$.

  Taking into account the equality
  $\frac{B(\nu-\mu,\mu)}{\nu}=\frac{B(\nu-\mu,\mu+1)}{\mu}$,
  from \eqref{61} and \eqref{61a}
  the equality \eqref{62} follows.

  The integral in the left hand side of \eqref{62a} is equal $I_1-I_2-I_3 (t^2+b)$,
  where
  \begin{equation*}
  \begin{split}
  &
  I_1=
  \nu\int_{t}^{+\infty}  S_{\nu}(y,u)  y|y^2-t^2|^{\nu-\mu-1}\;dy=
  \frac{B(\nu-\mu,\mu+1)}{2}\,\nu S_{\mu}(t,u)
  \,,
  \\&
  I_2=
  (\nu+1)\int_{t}^{+\infty}  S_{\nu+1}(y,u)  y|y^2-t^2|^{\nu-\mu}\;dy=
  \frac{B(\nu-\mu,\mu+1)}{2}\,(\nu-\mu) S_{\mu}(t,u)
  \,,
  \\&
  I_3=
  (\nu+1)\int_{t}^{+\infty}  S_{\nu+1}(y,u)  y|y^2-t^2|^{\nu-\mu-1}\;dy=
  \frac{B(\nu-\mu,\mu+1)}{2}\,(\mu+1) S_{\mu+1}(t,u)
  \,.
  \end{split}
  \end{equation*}
  To calculate $I_1$ we used \eqref{61a}.
  By changing  $\nu$ by $\nu+1$ in \eqref{61a} we calculated $I_2$.
  To calculate $I_3$, we changed  $\nu$ by $\nu+1$, and $\mu$  by $\mu+1$ in \eqref{61a}.
  The equality \eqref{62a} proved.
\end{proof}
For functions $h$ defined on $(0,+\infty )$, let us define a Hankel
transform\footnote{A non-standard notation is used here.} for
$m\in\C$ and $t>0$:
\begin{equation}\label{0c}
\F_m(h)(t):= t^{1-\frac{m}{2}}
             \int_{0}^{+\infty} h(x) x^{\frac{m}{2}} J_{\frac{m}{2}-1} (tx)\,dx=
             \int_{0}^{+\infty} h(x) x^{m-1} j_{\frac{m}{2}-1} (tx)\,dx
             \,,
\end{equation}
 where $J_{\lambda}$ is a Bessel function of the first kind, and
\begin{equation}\label{bes}
j_{\lambda}(x):=\frac{J_{\lambda}(x)}{x^{\lambda}}
 =\frac{1}{2^{\lambda}} \sum_{k=0}^{\infty}
 \frac{1}{\Gamma(k+\lambda+1)}\cdot
 \frac{\left(-\frac{x^2}{4}\right)^{k}}{k!}\;,\;x\in\C\;,\;\lambda\in\C\,.
\end{equation}
For $m\in\N$ the transform $\F_m$ is related to a Fourier transform
for radial functions:
 \begin{equation}\label{0f}
 F_m(h(||\cdot||_2))(x)=(2\pi)^{\frac{m}{2}}
 \F_m(h)(||x||_2) \,, \,x\in\R^m\,.
 \end{equation}
\begin{lemma}
 Let
\begin{equation}\label{66}
\begin{split}
  &
  g_{p,u}(x):=\frac{e^{-p\,x}}{x}-\frac{e^{-ux}}{e^x-1}\;\;,\;\;
  g_{p,u}(0):=u+\frac 12-p\;\;;\;\;\\&
  G_{p,u,\mu}(x):=xg'_{p,u}(x)+(2\mu-1)g_{p,u}(x)\;;
  \\&
  h_{u}(x):=\frac{e^{-ux}\, x}{e^x-1}\;\;,\;\;h_{u}(0):=1\,.
  \end{split}
\end{equation}
 Then equalities
\begin{equation}\label{67}
 \frac{1}{\mu(p^2+t^2)^{\mu}}-S_{\mu}(t,u)=
 \frac{\sqrt{\pi}\;\F_{2\mu+1}(g_{p,u})(t)}{2^{\mu-\frac{1}{2}}\,\Gamma(\mu+1)}\;
 \;,\;t\ge 0\,,\,p>0\,,\,u>-1\,,\,\mu>0\,,
\end{equation}
\begin{equation}\label{68}
 \frac{1}{\mu(p^2+t^2)^{\mu}}-S_{\mu}(t,u)=
 \frac{\sqrt{\pi}\;\F_{2\mu-1}(G_{p,u,\mu})(t)}{2^{\mu-\frac{1}{2}}\,\Gamma(\mu+1)\,t^2}\;
 \;,\;t>0\,,\,p>0\,,\,u>-1\,,\,\mu>\frac{1}{2}\,,
\end{equation}
\begin{equation}\label{69}
 \frac{d}{dt}\left\{t^{2\mu}S_{\mu}(t,u)\right\}=
 -\frac{\sqrt{\pi}\,t^{2\mu-1}\;\F_{2\mu+1}(h'_{u})(t)}{2^{\mu-\frac{1}{2}}\,\Gamma(\mu+1)}\;
 \;,\;t>0\,,\,u>-1\,,\,\mu>0
\end{equation}
hold.
\end{lemma}
\begin{proof}
  By \cite[Chapter 4.14(7,8), $\nu=\mu-\frac 12$]{Bateman_tables}
  we get the equalities
\begin{equation}\label{610}
 \frac{1}{\mu(p^2+a^2)^{\mu}}=
 \frac{\sqrt{\pi}}{2^{\mu-\frac{1}{2}}\,\Gamma(\mu+1)}
 \int_{0}^{+\infty}e^{-p\,x} x^{2\mu-1}j_{\mu-\frac 12}(ax)\,dx
 \;\;,\;a\in\R\,,\,p>0\,,\,\mu>0\,,
\end{equation}
\begin{equation}\label{611}
 \frac{2p}{(p^2+a^2)^{\mu+1}}=
 \frac{\sqrt{\pi}}{2^{\mu-\frac{1}{2}}\,\Gamma(\mu+1)}
 \int_{0}^{+\infty}e^{-p\,x} x^{2\mu}j_{\mu-\frac 12}(ax)\,dx
 \;\;,\;a\in\R\,,\,p>0\,,\,\mu>-\frac 12\,.
\end{equation}
 After the summation of the equalities  \eqref{611} for $p=k+u$, $k\in\N$, $u>-1$, $a=t\ge 0$,
 $\mu>0$, we get
\begin{equation}\label{612}
 S_{\mu}(t,u)=
 \frac{\sqrt{\pi}}{2^{\mu-\frac{1}{2}}\,\Gamma(\mu+1)}
 \int_{0}^{+\infty}\frac{e^{-ux}}{e^x-1}\, x^{2\mu}j_{\mu-\frac 12}(tx)\,dx
 \;\;,\;t\ge 0\,,\,u>-1\,,\,\mu>0\,.
\end{equation}
 For $a=t\ge 0$, the equalities \eqref{612} and \eqref{610} imply the equality \eqref{67}.
 If, in addition,  $\mu>\frac 12$, then for $t>0$ we can use integration in part in the right hand side of
 \eqref{67}. Taking into account
 $\frac{d}{dx}(j_{\lambda}(tx))=-t^2xj_{\lambda+1}(tx)$,
 we get
 \begin{equation*}
 \F_{2\mu+1}(g_{p,u})(t)=-\frac{1}{t^2}\int_{0}^{+\infty}g_{p,u}(x)x^{2\mu-1}\,dj_{\mu-\frac 32}(tx)
 =
 \frac{1}{t^2}\int_{0}^{+\infty}\left(g_{p,u}(x)x^{2\mu-1}\right)' j_{\mu-\frac 32}(tx)\,dx
 \,,
 \end{equation*}
 from which \eqref{68} follows.

 Let $h'_{u}(x):=\frac{d}{dx}(h_u(x))$. If $t>0$, by the substitution $s=tx$, in the integral from \eqref{612}, for $u>-1$,
 $\mu>0$,
 we get
\begin{equation*}
 \frac{d}{dt}\left\{t^{2\mu}S_{\mu}(t,u)\right\}=
 \frac{\sqrt{\pi}}{2^{\mu-\frac{1}{2}}\,\Gamma(\mu+1)}
 \int_{0}^{+\infty}h'_{u}\left(\frac{s}{t}\right)\, \left(-\frac{s}{t^2}\right)\,s^{2\mu-1}j_{\mu-\frac 12}(s)\,ds
 \,.
\end{equation*}
 The last equality is equivalent to \eqref{69}.
\end{proof}
\begin{lemma}\label{lem63}
 Let $p,u\in\R$. Then, for
  $ g_{p,u}(x):=\frac{e^{-p\,x}}{x}-\frac{e^{-ux}}{e^x-1}$, $x>0$,
  the next statements hold:
  \begin{enumerate}
  \item If $p-u\le\frac 12$ then $ g_{p,u}(x)>0$, for every $x>0$.
  \item If $p-u\ge 1$ then $ g_{p,u}(x)<0$, for every $x>0$.
  \item If $\frac 12< p-u< 1$ then, on an interval $(0,+\infty)$, the function $ g_{p,u}(x)$ is not of constant signs.
  \end{enumerate}
\end{lemma}
\begin{proof}
 From the equality $g_{p,u}(x)=\frac{e^{-ux}}{x(e^x-1)}\,f_{p-u}(x)$,
 where
 $f_a(x)=e^{-ax}(e^x - 1)-x$, follows that we have to find such a values of  $a\in\R$,
 that $f_a(x)$ preserves sign over the interval  $(0,+\infty)$.
 It is obvious that $f_a(0)=0$, $f'_a(x)=(1-a)e^{(1-a)x}+ae^{-ax}-1$,
 $f'_a(0)=0$, and  $f''_a(x)=((1-a)^2 e^{x}-a^2)e^{-ax}$.

 If $0\le a\le\frac 12$, $x>0$, then
 $e^x>1\ge\frac{a^2}{(1-a)^2}$. Consequently, $f''_a(x)>0$ and
 $f_a(x)>0$ for every $x>0$.
 If $a<0$,  $x>0$, then $f_a(x)>e^x-1-x>0$.

 If $a=1$, then $f''_1(x)=-e^{-x}<0$   and, consequently,
 $f_1(x)<0$ for every $x>0$.
 If $a>1$, then  $f_a(x)<f_1(x)<0$ for every $x>0$.

 If $\frac 12<a<1$, then $f_a(+\infty)=+\infty$ and $f''_a(0)=1-2a<0$.
 Consequently, $ f_a(x)$ does not preserve a sign on  $(0,+\infty)$.
\end{proof}
\begin{lemma}\label{le64}
  \begin{enumerate}
 \item
 If for some   $\nu>0$, $b\ge 0$, $u\in\R$, the function
 $(t^2+b)^{\nu}S_{\nu}(t,u)$ increases (decreases) with respect to   $t\in[0,+\infty)$,
 then,
 for every  $0<\mu\le\nu$ the function
 $(t^2+b)^{\mu}S_{\mu}(t,u)$ strictly increases (strictly decreases) with respect to  $t\in[0,+\infty)$.
 \item
 If for some  $\nu>0$, $b\ge 0$, $u\ge -1$, the function
 $(t^2+b)^{\nu}S_{\nu}(t,u)$ increases with respect to $t\in[0,+\infty)$,
 then,
 for every $0<\mu\le\nu$, $0\le a\le b$ the function
 $(t^2+a)^{\mu}S_{\mu}(t,u)$ strictly increases with respect to  $t\in[0,+\infty)$
 and, for $t\ge 0$, the inequality $\mu(t^2+a)^{\mu}S_{\mu}(t,u)<1$ take place.
 If, in addition,  $(\mu,a)\ne (\nu,b)$ then, for $t>0$, the
 inequality
  $\frac{d}{dt}\left((t^2+a)^{\mu}S_{\mu}(t,u)\right)>0$
  take place.
 \end{enumerate}
\end{lemma}
\begin{proof}
  The first statement we will prove by contradiction.Let us suppose that the function
   $(t^2+b)^{\nu}S_{\nu}(t,u)$ is monotonic with respect to $t\in[0,+\infty)$ but not
   strictly monotonic on this interval.
  Then
  \begin{equation}\label{Hnub}
  \begin{split}
  &
 \frac{d}{dt}\left((t^2+b)^{\nu}S_{\nu}(t,u)\right)=
 2t(t^2+b)^{\nu-1}H_{\nu,b}(t,u)\;,\;t>0\;,\;\text{ where }
 \\&
 H_{\nu,b}(t,u)=\nu S_{\nu}(t,u)-(\nu+1) S_{\nu+1}(t,u)(t^2+b)\;,\;t\in\R\;,
 \end{split}
  \end{equation}
  preserve sign for $t>0$ and vanishing on some interval.
  The function  $H_{\nu,b}(t,u)$ is analytic, for every fixed $u\in\R$, in the
  neighborhood of any real point, the, by the uniqueness theorem for analytic functions,
   $H_{\nu,b}(t,u)=0$ for every $t\in\R$.
   The asymptotic representation \eqref{asSmu} implies that
  \begin{equation}
 H_{\nu,b}(t,u)\sim
 \sum_{k=1}^{\infty}
 \frac{B_{2k}(-u)-b\,B_{2k-2}(-u)}{t^{2(k+\nu)}}\cdot
 \frac{(-1)^{k-1}\Gamma(\nu+k)}{\Gamma(\nu+1)\Gamma(k)}
 \;,\;t\to+\infty\;,
  \end{equation}
  and, consequently, $B_{2k}(-u)-b\,B_{2k-2}(-u)=0$ for every $k\in\N$.
  It is easy to check that the system of equations has no solutions
  (it is enough to check it for  $k=1,2,3$). Statement {\bf 1},
  for $\mu=\nu$ proved. If $0<\mu<\nu$, then the equality \eqref{62a} implies that
  the sign of $H_{\mu,b}(t,u)$ is equal to the sign of $H_{\nu,b}(t,u)$ for $t>0$.

  Let us prove the statement {\bf 2}. Let us suppose that the function $(t^2+b)^{\nu}S_{\nu}(t,u)$
  increasing with respect to $t\in[0,+\infty)$ for some  $\nu>0$, $b\ge 0$, $u\ge -1$.
  From what we already proved, this function increases  with respect to
  $t\in[0,+\infty)$ and, consequently, $H_{\nu,b}(t,u)\ge 0$ for $t>0$.
  The condition  $u\ge-1$ implies that $S_{\mu}(t,u)>0$ for every $\mu>0$, $t\in\R$.
  Consequently, the inequality $$H_{\nu,a}(t,u)=H_{\nu,b}(t,u)+(b-a)(\nu+1)S_{\nu+1}(t,u)\ge 0$$
  holds for every  $0\le a\le b$ and $t>0$.
  If, in addition, $a<b$ then $H_{\nu,a}(t,u)>0$ for $t>0$.
  If $0<\mu<\nu$, $0\le a\le b$ then the equality \eqref{62a} implies that
  $H_{\mu,a}(t,u)>0$ for $t>0$.

  Thus, for every $0<\mu\le\nu$, $0\le a\le b$ the function
 $(t^2+a)^{\mu}S_{\mu}(t,u)$ is strictly increasing with respect to $t\in[0,+\infty)$ and,
 consequently, the inequality
 $$\mu(t^2+a)^{\mu}S_{\mu}(t,u)<\mu\lim\limits_{x\to+\infty}(x^2+a)^{\mu}S_{\mu}(x,u)=1$$
 holds for $t\ge 0$.
  The statement {\bf 2} proved.
  \end{proof}
\begin{proof}[Proof of Theorem {\rm\ref{thSmu1}}]
  From Poisson's representation
   (see, for example, \cite[\S 7.3]{Bateman})
   \begin{equation*}
 j_{\lambda}(t)=\frac{1}{2^{\lambda}\Gamma\left(\lambda+\frac 12\right)\sqrt{\pi}}
 \int_{-1}^{1}e^{-itu}(1-u^2)^{\lambda-\frac 12}\,du\;,\;
 \lambda>-\frac 12\,,\,t\in\R
 \end{equation*}
 follows the inequality
 $|j_{\lambda}(t)|<j_{\lambda}(0)$, $\lambda>-\frac 12$, $t>0$.
 Further, for $t\ge 0$, $p>0$, $u>-1$, $\mu>0$,
 we use the equality \eqref{67} and the fact that, in the cases  $p-u\le\frac 12$ or $p-u\ge 1$,
 by the Lemma~\ref{lem63}, the function $g_{p,u}(x)$ preserves sign for $x>0$.
 In these cases, for every $t>0$, the inequality
 $|\F_{2\mu+1}(g_{p,u})(t)|<|\F_{2\mu+1}(g_{p,u})(0)|$ holds and is equivalent to the inequality \eqref{ner}.
 The Theorem~\ref{thSmu1} proved.
 \end{proof}
\begin{proof}[Proof of Theorem  {\rm\ref{thSmu2}}]
 The relations \eqref{asSmu} and inequality \eqref{2} imply that, for every fixed $\mu>0$ and $u\ge 0$,
 the function $S(t)=S_{\mu}(t,u)$, $t\ge 0$,
 satisfies to the conditions of the Theorem \ref{thSt} for $C=\frac{1}{\mu}$,
 $\delta_1=2\mu$, $A=u^2+u+\frac 16$, $\beta_1=2$.
 From the Theorem \ref{thSt} (statement~{\bf 3}) follows that
 $f_{\mu,u}(+\infty)=u^2+u+\frac{1}{6}$ and, for mentioned  $a$ and $b$
 the inequality \eqref{ner2} take place.

  Let us prove the statement {\bf 1.} The inequalities $u^2+u<m_{\mu}(u)\le M_{\mu}(u)<(1+u)^2$,
  for $\mu>0$, $u\ge 0$, follows from the inequality \eqref{mM1} and Remark \ref{re2}
  (the cases~1 and~2) for the function  $g(x)=x^{-\mu-1}$, taking into account that
   $S_{\mu}(t,u)=S(u,t^2)$,  $f_{\mu,u}(t)=\psi(u,t^2)$,
  $m_{\mu}(u)=m(u)$ и $M_{\mu}(u)=M(u)$ (see~\eqref{mM2}).\\
  To prove the inequality $u^2+u+\frac{1}{4}<M_{\mu}(u)$ we use the inequality \eqref{67}
  and the Theorem from~\cite{Zast2004}.
  Let us define, for  $m\in\R$,
  \begin{equation}
  \begin{split}
  &
  H_m:=L\left((0,1),x^{m-1}dx\right)\cap
 L\left((1,+\infty),x^{\frac{m-1}{2}}dx\right)\;,
 \\&
 H_m^1:=L\left((0,+\infty),x^{m-1}dx\right)\;,\;
 H_m^2:=L\left((0,+\infty),x^{\frac{m-1}{2}}dx\right)\;.
\end{split}
 \end{equation}
\begin{theorem}[Zastavnyi \cite{Zast2004}]\label{thZast2004}
 \begin{enumerate}
 \item  Let $m>0$ and $f\in \left(H_{m}^1\cup H_{m}^2\right)\cap H_{m}$.
 If, for some fixed $\delta>0$ and $\varepsilon\in\C$, the function $f$
 is bounded on $(0,\delta)$ and $|\F_m(f)(t)|-\varepsilon\F_m(f)(t)\in
 H_{m}^1$, then $\F_m(f)\in H_{m}^1$.
 \item  Let $m>0$ and $f\in \left(H_{m}^1\cup H_{m}^2\right)\cap
H_{m}$. If $\F_m(f)(t)$ preserves sign for $t>0$ and $f$ is
continuous at zero then $\F_m(f)\in H_{m}^1$. If, in addition,
$f(0)=0$ then $f(x)=0$ for almost all $x>0$.
 \end{enumerate}
\end{theorem}
 Let $\mu>0$, $u\ge 0$, and $p=u+\frac 12$. Let us suppose that
$\F_{2\mu+1}(g_{p,u})(t)$ preserves sign for $t>0$ where the
function $g_{p,u}$ defines by equality~\eqref{66}. Then the function
$f(x):=g_{p,u}(x)$, $x\ge 0$, satisfies to the conditions of
Theorem~\ref{thZast2004} (for $m=2\mu+1$, $f\in H_m^1\cap H_{m}$).
Moreover, $f$ is continuous at zero and $f(0)=0$. Therefore,
$f(x)=0$ for almost all  $x>0$. But it is impossible.
 Therefore (see equality~\eqref{67}) the difference
 $\frac{1}{\mu((u+\frac 12)^2+t^2)^{\mu}}-S_{\mu}(t,u)$
 does not preserve sign for $t>0$ and, consequently, $u^2+u+\frac{1}{4}<M_{\mu}(u)$.
 The statement {\bf 1} of Theorem~\ref{thSmu2} proved.

 The statement {\bf 2} of Theorem~\ref{thSmu2} immediately follows from the equality~\eqref{62}.

 Let us prove the statement {\bf 3.} Let $\nu>0$, $u\ge 0$. Then
 \begin{equation*}
 f'_{\nu,u}(t)=2t
 \left(
 \left(\nu S_{\nu}(t,u)\right)^{-\frac{1}{\nu}-1}
 (\nu+1)S_{\nu+1}(t,u)-1
 \right)
 \;.
 \end{equation*}
 The last equality imply that for every $0<\alpha<\beta$
 the next three conditions
 \begin{subequations}\label{ekv}
 \begin{equation}
 f'_{\nu,u}(t)\le(<)0\;,\;\alpha<t<\beta\;,
 \end{equation}
 \begin{equation}
 \psi_{\nu,u}(t):=
  \left((\nu+1)S_{\nu+1}(t,u)\right)^{\frac{1}{\nu+1}}-
   \left(\nu S_{\nu}(t,u)\right)^{\frac{1}{\nu}}
 \le(<)0\;,\;\alpha<t<\beta\;,
 \end{equation}
 \begin{equation}
 f_{\nu,u}(t)-f_{\nu+1,u}(t)\le(<)0\;,\;\alpha<t<\beta
 \end{equation}
 \end{subequations}
 are equivalent.
 From the expansion
$(1+x)^{\frac{1}{\mu}}=1+\frac{x}{\mu}+\frac{x^2(1-\mu)}{2\mu^2}+O(x^3)$,
$x\to 0$, and asymptotic~\eqref{asSmu} follows that the expansion
\begin{equation*}
 \left(\mu S_{\mu}(t,u)\right)^{\frac{1}{\mu}}
 =\frac{1}{t^2}-\frac{B_2(-u)}{t^4}+\frac{B_4(-u)(\mu+1)+B_2^2(-u)(1-\mu)}{2t^6}+O\left(\frac{1}{t^8}\right)
 \;,\;t\to+\infty
\end{equation*}
 takes place for every $\mu>0$. Therefore,
 \begin{equation}
 \psi_{\nu,u}(t)=\frac{B_4(-u)-B_2^2(-u)}{2t^6}+O\left(\frac{1}{t^8}\right)\sim
 -\,\frac{60u^2+60u+11}{360t^6}
 \;,\;t\to+\infty\;.
 \end{equation}
  This asymptotic implies the existence such a value  $t_{\nu,u}\ge 0$ that,
  for every $t>t_{\nu,u}$, the inequality  $\psi_{\nu,u}(t)<0$ holds and, consequently, the inequality
  $f'_{\nu,u}(t)<0$ holds.
  Therefore, the function $f_{\nu,u}(t)$ strictly decreases with
  respect to $t\in[t_{\nu,u},+\infty)$.
  Consequently, for every $t>t_{\nu,u}$, the inequality
  $f_{\nu,u}(+\infty)<f_{\nu,u}(t)<f_{\nu,u}(t_{\nu,u})$ holds which is equivalent to the inequality~\eqref{ner1}.
  The statement {\bf 3} of Theorem~\ref{thSmu2}  proved.

  Let us prove statement {\bf 4.} Let, for some $\nu>0$, $u\ge 0$, the inequality \eqref{ner3} be satisfied.
  From the equivalence of the of the conditions~\eqref{ekv},
  we get that, for every $t>0$, the inequality $f'_{\nu,u}(t)\le 0$
  be satisfied and, hence, the function $f_{\nu,u}(t)$ degreases on
  $[0,+\infty)$.
  Since  $f_{\nu,u}(t)$ is analytic function in the neighborhood  of every real point then
  $f_{\nu,u}(t)$ strictly decreases with respect to  $t\in[0,+\infty)$
  (in the opposite $f'_{\nu,u}(t)= 0$ on some interval and, consequently, for every  $t\in\R$
  which contradict to the statement~{\bf 3}).
  Therefore, the inequality~{\rm\eqref{ner1}} be satisfied for
   $t_{\nu,u}=0$, $M_{\nu}(u)=f_{\nu,u}(0)$, and
   $m_{\nu}(u)=f_{\nu,u}(+\infty)=u^2+u+\frac{1}{6}$.
  From the equivalence of the conditions~\eqref{ekv},
  we get, for every  $t>0$, the inequality
  $f_{\nu,u}(t)\le  f_{\nu+1,u}(t)$ be satisfied and, hence,
  $m_{\nu}(u)\le m_{\nu+1}(u)\le f_{\nu+1,u}(+\infty)=u^2+u+\frac{1}{6}$.
  Therefore,  $m_{\nu+1}(u)=u^2+u+\frac{1}{6}$.
  If $0<\mu< \nu+1$ then the statements~{\bf 1} and {\bf 2} imply the inequality
  $m_{\nu+1}(u)\le m_{\mu}(u)\le u^2+u+\frac{1}{6}$ and, hence,
  $m_{\mu}(u)= u^2+u+\frac{1}{6}$.

  Let us prove that for   $\mu\in(0, \nu+1]$ the inequality
  $S_{\mu}(t,u)<\frac{1}{\mu(t^2+m_{\mu}(u))^{\mu}}$ be satisfied for $t\ge 0$.
  For  $\mu=\nu+1$, this inequality follows from the inequality
  $f_{\nu,u}(t)\le f_{\nu+1,u}(t)$, $t\ge 0$, which we proved.
  We used that  $m_{\nu+1}(u)=f_{\nu,u}(+\infty)<f_{\nu,u}(t)$, $t\ge 0$.
  For $0<\mu<\nu+1$, we used the equality~\eqref{62} in which
  $p=m_{\mu}(u)= u^2+u+\frac{1}{6}$ and instead $\nu$ we put $\nu+1$.

  Let us prove that for every $0<\mu\le\nu$, $0\le a\le u^2+u+\frac{1}{6}$, and $t>0$,
  the inequality
  $\frac{d}{dt}\left((t^2+a)^{\mu}S_{\mu}(t,u)\right)>0$ be satisfied.
   It is already proved that the inequality~{\rm\eqref{ner1}} be satisfied for
   $t_{\nu,u}=0$. Therefore, for every $t>0$, the inequality
   $\nu S_{\nu}(t,u)<(t^2+b)^{-\nu}$, where $b=u^2+u+\frac{1}{6}$, be satisfied.
   Taking it and inequality~\eqref{ner3} into account, we get that
   for $t>0$ the inequality
   \begin{equation*}
   (\nu+1)S_{\nu+1}(t,u)\le
   \left(\nu S_{\nu}(t,u)\right)^{1+\frac{1}{\nu}}<
   \frac{\nu S_{\nu}(t,u)}{t^2+b}
   \end{equation*}
   be satisfied.
   The last inequality is equivalent (see~\eqref{Hnub}) to the inequality
   $$\frac{d}{dt}\left((t^2+b)^{\nu}S_{\nu}(t,u)\right)>0\;, \,t>0.$$
   Further, we apply the Lemma~\ref{le64}.
   The statement~{\bf 4} of Theorem~\ref{thSmu2} proved.

 The statement {\bf 5} immediately follows from the statement~{\bf 4},
 because the inequality~\eqref{ner3} proved in 1998, for $\nu=1$, $u=0$,
 by Wilkins~\cite{Wilkins}.

 The Theorem~\ref{thSmu2} proved.
\end{proof}
\begin{proof}[Proof of the Theorem {\rm\ref{thSmu3a}}]
   Let us first prove the statement {\bf 1.} It is easy to see that for  $p\ge 0$,
   $\mu>0$, the next equality
   \begin{equation}\label{613}
   (-1)^k\psi^{(k)}_{p,u,\mu}(t)=\frac{\Gamma(\mu+k+1)}{\Gamma(\mu+1)}\psi_{p,u,\mu+k}(t)
   \;,\;k\in\Z_+\;,\;t>0
   \end{equation}
   holds.
  If, for some $\mu>0$, $\psi_{p,u,\mu}\in M(0,+\infty)$ then $p\ge 0$
  and it follows from~\eqref{613} that the inequality $\psi_{p,u,\mu+k}(t)\ge 0$ holds for every
  $k\in\Z_+$ and $t>0$. Each of these inequalities is equivalent to the inequality
  $p\le m_{\mu+k}(u)$.
  Passing to the limit in the last inequality as $k\to+\infty$,
  we get  $p\le m_{\infty}(u)$.\\
  If $0\le p\le m_{\infty}(u)$ then $\psi_{p,u,\mu}(t)\ge 0$ for $\mu>0$, $t>0$ and consequently,
  by the equality~\eqref{613},
  $\psi_{p,u,\mu}\in M(0,+\infty)$ for $\mu>0$.
  The statement {\bf 1} proved.

  In an analogous way we prove the statement {\bf 2.}

  Let us prove the statement {\bf 3.} It is easy to check that for $p\ge 0$,
  $u\ge 0$,   $\mu>0$ the next equality
  \begin{equation*}
  \psi_{p,u,\mu}(t)=\frac{1}{\Gamma(\mu+1)}\int_{0}^{+\infty}e^{-tx}x^{\mu-1}(e^{-px}-\varphi_{u}(x))\,dx
  \;,\;t>0
  \end{equation*}
  holds.
  From this inequality and the Theorem Hausdorff-Bernstein-Widder follows that,
  for $p\ge 0$,  $u\ge 0$,   $\mu>0$, the condition  $\psi_{p,u,\mu}\in M(0,+\infty)$
   is equivalent to the inequality
   $e^{-px}-\varphi_{u}(x)\ge 0$, $x>0$. The last inequality is equivalent to
   $p\le -\sup_{x>0}\frac{\ln\varphi_{u}(x)}{x}$.
   Taking into the statement~{\bf 1}, the last inequality implies the equality
   $m_{\infty}(u)=-\sup_{x>0}\frac{\ln\varphi_{u}(x)}{x}$.
   By the same way, we can prove that
   $M_{\infty}(u)=-\inf_{x>0}\frac{\ln\varphi_{u}(x)}{x}$.

      Let us prove the statement {\bf 4.}
     The equality~\eqref{228} imply that,
     for  $u\ge 0$ and $x>0$, the strict inequalities
     $u^2+u<-\frac{\ln\varphi_{u}(x)}{x}<(u+1)^2$ holds.
     It is obviously that, for every fixed $u\ge 0$, the relation
     $\varphi_u(x)=2x(u+1)e^{-x(u+1)^2}(1+o(1))$, $x\to+\infty$,
     holds. This relation and the Example~\ref{ex1} imply that next equalities
     \begin{equation*}
     \lim_{x\to+0}-\frac{\ln\varphi_{u}(x)}{x}=u^2+u+\frac{1}{6}\;,\;
     \lim_{x\to+\infty}-\frac{\ln\varphi_{u}(x)}{x}=(u+1)^2
     \end{equation*}
     hold.
     Therefore,
     \begin{equation*}
     \begin{split}
     &
     m_{\infty}(u)=\inf_{x>0}-\frac{\ln\varphi_{u}(x)}{x}>u^2+u\;,\;
     m_{\infty}(u)\le u^2+u+\frac{1}{6}\;,
     \\&
     (u+1)^2=\lim_{x\to+\infty}-\frac{\ln\varphi_{u}(x)}{x}\le
     M_{\infty}(u)=\sup_{x>0}-\frac{\ln\varphi_{u}(x)}{x}\le (u+1)^2\;.
     \end{split}
     \end{equation*}
    The Theorem~\ref{thSmu3a} proved.
\end{proof}
\begin{proof}[Proof of the Theorem {\rm\ref{thSmu3}}]
 Let us prove the sufficiency in the statement {\bf 1.}

  Let $u\ge 0$ and $0\le p\le \sqrt{m_{\infty}(u)}$.
  From the Theorem~\ref{thSmu3a} follows that, $m_{\infty}(u)>0$.
  The equality~\eqref{67} and right hand side of ~\eqref{ner2} imply that,
  for every $\mu>0$ and $0< p\le \sqrt{m_{\infty}(u)}$,
  the inequality $\F_{2\mu+1}(g_{p,u})(t)\ge 0$, $t\ge 0$, holds.
  Next, we use the connection~\eqref{0f} between Hankel's transform
  and Fourier's transform of radial functions.
  Then, by the Bochner's Theorem, we obtain
  $g_{p,u}\in\Phi(l^{m+1}_{2})$, for each $m\in\N$ and $0<p\le\sqrt{m_{\infty}(u)}$.
  The convergence $g_{p,u}\to g_{0,u}$, as $p\to+0$, implies  $g_{0,u}\in\Phi(l^{m+1}_{2})$,
  for every $m\in\N$.
  Then, by  Schoenberg's Theorem, $g_{p,u}\in\Phi(l_{2})$ for every $0\le p\le\sqrt{m_{\infty}(u)}$.
  The sufficiency in the statement {\bf 1} proved.

     Let us prove the necessity in the statement {\bf 1}.
     Let $g_{p,u}\in\Phi(l_{2})$.
     Then, by Schoenberg's Theorem, $g_{p,u}\in\Phi(l^{m+1}_{2})$, for every $m\in\N$.
     Boundedness of every positive definite function implies $p\ge 0$.
     If $p>0$ then, by  Bochner's Theorem and equality~\eqref{0f}, the inequality
     $\F_{2\mu_m+1}(g_{p,u})(t)\ge 0$, where $t\ge 0$, holds, for every $\mu_m=\frac m2$, $m\in\N$.
     The equality~\eqref{67} and right hand side of~\eqref{ner2} imply
     $p\le\sqrt{m_{\mu_m}(u)}$, $m\in\N$.
     Passing to the limit, as $m\to\infty$, and taking into account~\eqref{235},
     we get $p\le\sqrt{m_{\infty}(u)}$.
     The necessity of statement {\bf 1} proved.

   We can prove the statement {\bf 2} in the same way that the statement {\bf 1} does.
   We consider the left hand side inequality in~\eqref{ner2} instead the right hand side in~\eqref{ner2}.

   The statement {\bf 3} follows from the next equivalent statements:\\
   {\bf a)} The inequality $\frac{d}{dt}\left\{t^{2\mu}S_{\mu}(t,u)\right\}\ge 0$
    holds, for every $t>0$ and $\mu>0$. \\
   {\bf b)} The inequality $\frac{d}{dt}\left\{t^{2\mu}S_{\mu}(t,u)\right\}\ge 0$
    holds, for every $t>0$ and $\mu=\frac m2$, $m\in\N$. \\
   {\bf c)} The inequality $-\F_{2\mu+1}(h'_{u})(t)\ge 0$ holds, for every
    $t>0$ and $\mu=\frac m2$, $m\in\N$.\\
   {\bf d)} $-\frac{d}{dx}(h_u(x))\in\Phi(l^{m+1}_{2})$, $m\in\N$.\\
   {\bf e)} $-\frac{d}{dx}(h_u(x))\in\Phi(l_{2})$.\\
   The Lemma~\ref{le64} implies the equivalence the statements {\bf a} and  {\bf b}.
   The equality~\eqref{69} implies the equivalence of  {\bf b} and  {\bf c}.
   The equality~~\eqref{0f} and Bochner's Theorem implies the equivalence of {\bf c} and  {\bf d}.
   The Schoenberg's Theorem implies the equivalence of {\bf d} and  {\bf e}.

 The Theorem~\ref{thSmu3} proved.
\end{proof}


\begin{thebibliography}{9}
\setlength{\itemsep}{5pt}
\bibitem{Mathieu}
  \'{E}.L. Mathieu,
    {\em Trait\'{e} de Physique Math\'{e}matique. VI-VII: Th\'{e}ory
                   de l’\'{E}lasticit\'{e} des Corps Solides (Part~2),}
                    Gauthier-Villars, Paris, 1890.

 \bibitem{Berg}
 L. Berg,
 ``{\"Uber eine Absch\"atzung von Mathieu}'',
 {\em Math. Nachr.\/},~{\bf 7} (1952), 257--259.

\bibitem{Schroder}
 K. Schr\"oder,
 ``{Das Problem der eingespannten rechteckigen elastishen Platte}'',
 {\em Math. Ann.\/},~{\bf 121} (1949), 247--326.

 \bibitem{Corput}
 J.G. {van der} Corput,  L.O. Heflinger,
 ``{On the inequality of Mathieu}'',
 {\em Indagationes Mathematicae\/},~{\bf 18} (1956), 15--20.

\bibitem{Emersleben}
 O. Emersleben,
 ``{\"Uber die Reihe $\sum k(k^2+c^2)^{-2}$}'',
 {\em Math. Ann.\/},~{\bf 125} (1952), 165--171.

\bibitem{Polya}
 G. Polya,
 ``{\"{U}ber die Nullstellen gewisser ganzer Funktionen}'',
 {\em Math. Z.\/},~{\bf 12} (1918), 352--383.

\bibitem{Lukacs}
E. Lukacs,
 {\em Characteristic functions},
 Griffin, London, 2nd edition, 1970.

\bibitem{Zast2003}
V.P. Zastavnyi,
 ``{Extension of a Function from the Exterior of an
   Interval to a Positive-Definite Function on the Entire Axis and an
   Approximation Characteristic of the Class $W_M^{r,\beta}$}'',
    {\em Ukrainian Mathematical Journal\/},~{\bf 55}(7) (2003),  1189--1197.

 \bibitem{Makai}
E. Makai,
 ``{On the inequality of Mathieu}'',
 {\em Publ. Math. Debrecen\/},~{\bf 5} (1957), 204-205.

 \bibitem{Elbert}
   A. Elbert,
  ``{Asymptotic expansion and continued fraction for Mathieu’s   series}'',
                    {\em Period. Math. Hungar.\/},~{\bf 13} (1982), 1–8.

\bibitem{Alzer}
 H. Alzer, J.L. Brenner, O.G. Ruehr,
 ``{On Mathieu's inequality}'',
 {\em J. Math. Anal. Appl.\/},~{\bf 218} (1998), 607-610.

\bibitem{Alzer97}
 H. Alzer, J.L. Brenner,
 ``{Problem 97-1}'',
 {\em SIAM Rev.\/},~{\bf 39} (1997), 123.

\bibitem{Wilkins}
 J.E. Wilkins, Jr,
 ``{An Inequality}'',
 {\em SIAM Rev.\/},~{\bf 40}(1) (1998), 126--128.

  \bibitem{Wang}
     C.L. Wang, X.H. Wang,
    ``{A refinement of the Mathieu inequality}'',
                {\em Univ. Beograd. Publ. Electrotehn. Fak. Ser. Mat. Fiz.}, No {\bf 716-734} (1981) 22-24.

\bibitem{Diananda}  P.H. Diananda,
  ``{Some Inequalities Related to an Inequality of Mathieu}'',
{\em Math. Ann.\/},~{\bf 250} (1980), 95--98.

 \bibitem{Hoorfar and Qi}   A. Hoorfar, F. Qi,
     ``{Some new bounds for Mathieu's series}'',
     {\em Abstract and Applied Analysis\/},~{\bf 2007} (2007), Article ID 94854,
      10  pages doi:10.1155/2007/94854.

 \bibitem{Qi}
  F. Qi,  CH.-P. Chen,  B.-N. Guo,
   ``{Notes on double inequalities of Mathieu’s series}'',
                {\em International Journal of Mathematics and Mathematical Sciences\/},~{\bf 16} (2005), 2547–2554.

 \bibitem{Tomovski}
  \v{Z}. Tomovski,
  ``{New double inequalities for Mathieu type series}'',
                     {\em Univ. Beograd. Publ. Elektrotehn. Fak. Ser. Mat.\/},~{\bf 15} (2004), 79–83.

 \bibitem{Bateman}
        H. Bateman, A. Erd\'{e}lyi,
                 {\em Higher transcendental functions, vol. 1, 2,}
                                   New York, Toronto, London, MC Graw-Hill Book Company, INC, 1953.

 \bibitem{Pogany}
  T.K. Pog\'{a}ny,  H.M. Srivastava,  \v{Z}. Tomovski,
 ``{Some families of  Mathieu {\bf a}-series and alternating Mathieu {\bf a}-series}'',
                  {\em Applied Mathematics and Computation\/},~{\bf 173} (2006), 69–108.

\bibitem{Zast2000} V.P. Zastavnyi,
  ``{On positive definiteness of some functions}'',
             {\em Journal of Multivariate Analysis\/},~{\bf 73} (2000), 55-81.

\bibitem{TrBel} R. M. Trigub,  E. S.  Belinsky,
 {\em Fourier Analysis and Approximation of Functions\/}, Boston,
 Dordrecht, London,
Kluwer-Springer, 2004.

\bibitem{Vakh}  N.N. Vakhania, W.I. Tareladze, S.A. Chobanjan,
{\em Probability distributions on Banach Spaces\/}, Moscow, Nauka,
1985 (in Russian).

\bibitem{Fell} W. Feller,
 ``{Completely monotone functions and sequences}'',
{\em  Trans. Amer. Math. Soc.\/},~{\bf 5} (1939), 662-774.

\bibitem{Scho1} I. J. Schoenberg,
  ``{Metric spaces and completely monotone functions}'',
 {\em  Ann. Math.\/},~{\bf 39} (1938), 811-841.

\bibitem{Scho2} I. J. Schoenberg,
   ``{Metric spaces and positive definite functions}'',
 {\em  Trans. Amer. Math. Soc.\/},~{\bf 44} (1938), 522-536.

\bibitem{Fedoryuk}  M. V. Fedoryuk,
  {\em Asymptotics: Integrals and series,\/}  Moscow, Nauka, 1987 (in Russian).

\bibitem{Bruijn} N.G. {de} Bruijn,
  {\em Asymptotic methods in analysis,\/}
  North-Holland Publishing Co., Amsterdam, 1958.

 \bibitem{Bateman_tables}
        H. Bateman, A. Erd\'{e}lyi, W. Magnus, F. Oberhettinger, F. Tricomi,
                 {\em Tables of integral transforms, vol. I,}
                                   New York, Toronto, London, MC Graw-Hill Book Company, INC, 1954.

\bibitem{Zast2004} V.P. Zastavnyi
                 ``{Sufficient
conditions of an integrability of the Hankel transformation}'',
  {\em Transactions of IAMM of NAS of Ukraine\/},~{\bf 9} (2004), 68-75 (in Russian).


\end{thebibliography}
\end{document}